\documentclass[11pt]{article}

\usepackage{bm}%
\usepackage[colorlinks=true,linkcolor=blue]{hyperref}%
\usepackage[T1]{fontenc}
\usepackage[utf8]{inputenc}
\usepackage{algorithm}
\usepackage[font={small,it}]{caption}
\usepackage{authblk}
\usepackage{float}
\usepackage{mathtools}
\usepackage{amsmath}
\usepackage{amsthm}
\usepackage{amssymb}
\usepackage{graphicx}
\usepackage{algpseudocode}
\makeatother
\usepackage{subfig}

\providecommand{\tabularnewline}{\\}
\floatstyle{ruled}
\newfloat{algorithm}{tbp}{loa}
\providecommand{\algorithmname}{Algorithm}
\floatname{algorithm}{\protect\algorithmname}

\theoremstyle{definition}
\newtheorem{defn}{\protect\definitionname}
\theoremstyle{plain}
\newtheorem{thm}{\protect\theoremname}
\theoremstyle{remark}
\newtheorem{rem}{\protect\remarkname}

\providecommand{\definitionname}{Definition}

\providecommand{\theoremname}{Theorem}

\makeatother

\providecommand{\definitionname}{Definition}
\providecommand{\remarkname}{Remark}
\providecommand{\theoremname}{Theorem}

\expandafter\ifx\csname package@font\endcsname\relax\else
 \expandafter\expandafter
 \expandafter\usepackage
 \expandafter\expandafter
 \expandafter{\csname package@font\endcsname}%
\fi
\author{Mattia Serra\thanks{serram@ethz.ch}}
\author{George Haller\thanks{Email address for correspondence: georgehaller@ethz.ch}}
\title{Objective Eulerian Coherent Structures}%
\affil{Institute for Mechanical Systems, ETH Zürich, Leonhardstrasse 21, 8092 Zurich, Switzerland}%
\begin{document}

\maketitle

\date{}

\begin{abstract}
	We define objective Eulerian Coherent Structures (OECSs) in two-dimensional,
	non-autonomous dynamical systems as instantaneously most influential
	material curves. Specifically, OECSs are stationary curves of the
	averaged instantaneous material stretching-rate or material shearing-rate
	functionals. From these objective (frame-invariant) variational principles,
	we obtain explicit differential equations for hyperbolic, elliptic
	and parabolic OECSs. As illustration, we compute OECSs in an unsteady
	ocean velocity data set. In comparison to structures suggested by
	other common Eulerian diagnostic tools, we find OECSs to be the correct
	short-term cores of observed trajectory deformation patterns.
\end{abstract}

\section{Introduction}

Coherent structures behind material pattern formation in unsteady
flows have received considerable attention (see, e.g., \cite{ProvenzaleBarotropicVortices1999,Boffetta_transpBarr2001,PeacockLCSFocus2010},
for reviews). Among other mathematical advances (\cite{Froyland2014,Allshouse2012,Budisic2012,Ma2014,Peacock2015}),
variational techniques have been developed to identify Lagrangian
Coherent Structures (LCSs), the most influential cores of material
deformation, in experimental and numerical velocity data (\cite{Haller2011VARIHyp,Farazmand2012,Farazmand2012a,BlackHoleHaller2013,ShearlessBarrierFarazmand2014}).

Specifically, initial positions of hyperbolic LCSs (generalized stable
and unstable manifolds), elliptic LCSs (generalized KAM tori) and
parabolic LCSs (generalized jet cores) are now readily computable
as solutions of differential equations defined over the initial configuration
of a system \cite{LCSHallerAnnRev2015}. Later positions of these
LCSs are then obtained by advecting their initial positions under
the flow map. By the objectivity (frame-invariance) of their underlying
variational principles, variational LCSs transform properly under
coordinate changes $x{\mapsto}{\tilde x}$ of the form 
\begin{equation}
x=Q(t)\tilde{x}+b(t),\label{eq:observer change}
\end{equation}
where $Q(t)$ is an arbitrary proper orthogonal tensor family generating
time-dependent rotations, and $b(t)$ is an arbitrary vector family
introducing time-dependent translations of the frame \cite{TruesdellNoll2004}. 

All LCSs, however, are intrinsically tied to a specific finite time
interval over which they exert their influence on nearby trajectories.
This influence is an integrated one, filtering out short-term anomalies
in the flow. Yet short-term variability in material structures is
often seen as significant in highly unsteady flows, explaining the
popularity of Eulerian (i.e., instantaneous velocity-based) diagnostic
fields in fluid dynamics (see, e.g., \cite{Jeong&Uss.On_Id_of_Vtx_lambda2_1995,HallerObjDefVtx2005},
and \cite{Chakraborty2005} for reviews).

Eulerian diagnostics generally highlight features of the velocity
field in a given frame of reference. Most studies of flows, nevertheless,
would ideally want to understand, forecast or control the evolution
and mixing of \emph{material trajectories}, as well as the transport
of the physical properties they carry. Indeed, the motivation for
each commonly used Eulerian diagnostic is invariably grounded in the
desire to understand particle motion over short time scales. Yet the
approximations and heuristics employed in deriving these diagnostics,
as well as the expectation for simply implementable results, tends
to change the original focus of these approaches from material features
to the analysis of frame-dependent velocity features. 

Even short-term identifications of \emph{material} coherence, however,
must be frame-invariant by definition. This is because material coherent
structures are composed of trajectories, which do not depend on coordinates
and hence must transform properly under \eqref{eq:observer change}
from the frame of one observer to the other's. If a proposed coherent
structure criterion labels different material sets as coherent in
different frames, then it either lacks a proper physical foundation,
or its physical foundation is implemented through erroneous mathematics.
Either way, the criterion cannot be used reliably in now-casting or
real-time decision-making under unsteady flow conditions \cite{Peacock2015}. 

The few available objective Eulerian coherent structure diagnostics
include that of Tabor and Klapper \cite{TaborKlapper1994}, who consider
a point elliptic (i.e, part of a vortex) if the vorticity expressed
in the basis of the rate-of-strain eigenvectors dominates strain (see
also \cite{Lapeyre1999,Lapeyre2001}). Another objective Eulerian
technique is the recent global approach \cite{HallerLAVD2015} to
rotationally coherent vortices. Identified as outermost closed level
sets of the instantaneous vorticity deviation (IVD) from the spatial
mean vorticity, IVD-based vortex boundaries are objective, short-term
limits of closed curves obtained from a global Lagrangian rotational
coherence principle. 

Here we introduce a more general approach to Objective Eulerian Coherent
Structures (OECSs) in two dimensions. We effectively seek these structures
as short-term limits of LCSs. In contrast to \cite{HallerLAVD2015},
we use stretching-based variational LCS theories in taking this limit.
As a consequence, we obtain a broader class of OECSs that includes
elliptic (vortex-type), hyperbolic (stretching or contracting) and
parabolic (jet-type) structures. 

We give an explicit parametrization of all these types of OECSs as
solutions of autonomous ordinary differential equations (ODEs). Our
approach is global, eliminating the shortcomings of local structure
identification schemes noted in \cite{LugtDilemmaVortxDef1979}. All
OECSs are computable without dependence on any chosen time scale.
They are also objective, depending solely on invariants of the rate-of-strain
tensor, i.e., the symmetric part of the velocity gradient. 

\begin{figure}[H]
	\subfloat[]{\includegraphics[width=0.47\textwidth]{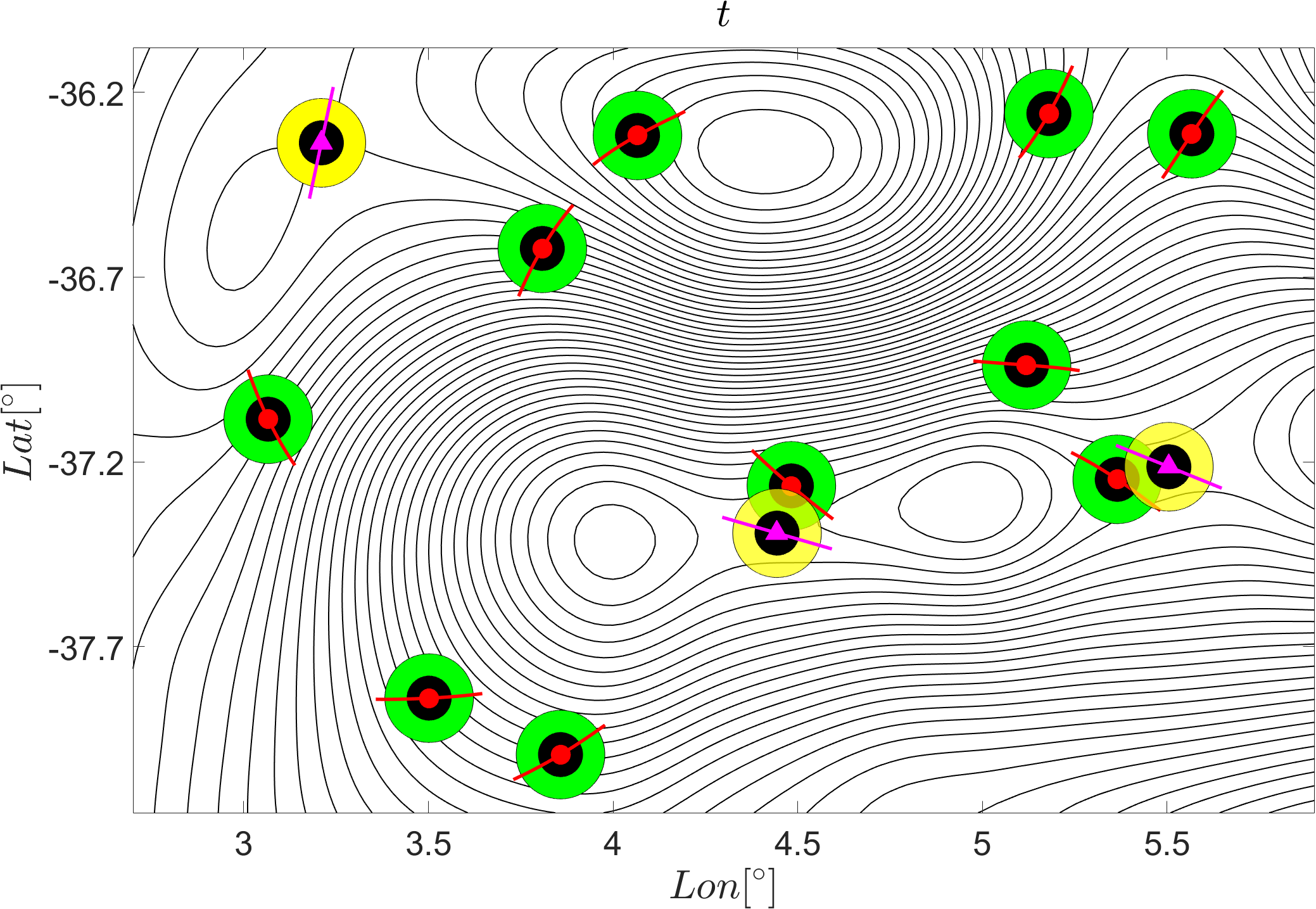}\label{fig:HypOECSsMotOrFrt0}}\hfill{}
	\subfloat[]{\includegraphics[width=0.47\textwidth]{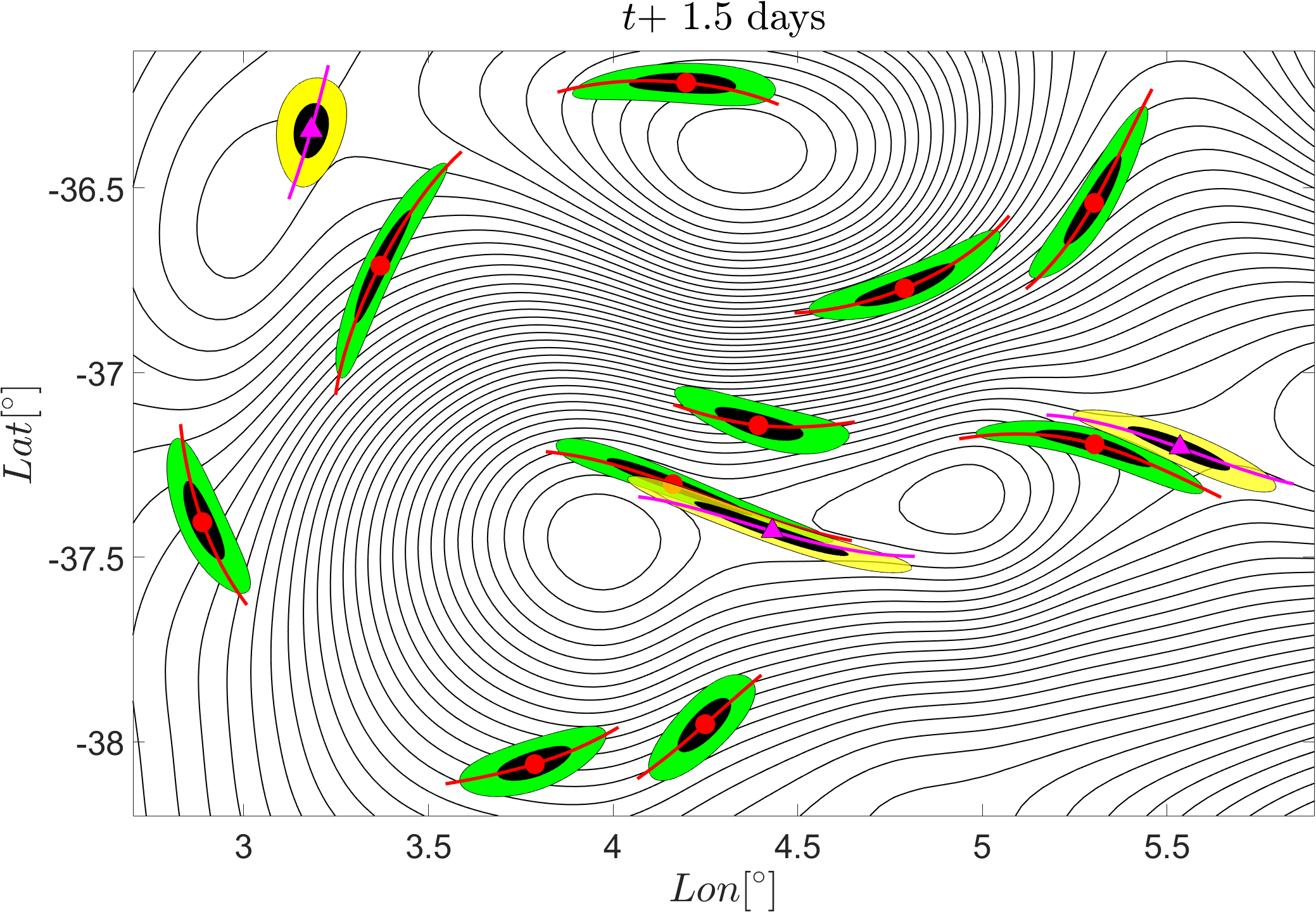}\label{fig:HypOECSsMotOrFrtf}}\\
	\subfloat[]{\includegraphics[width=0.47\textwidth]{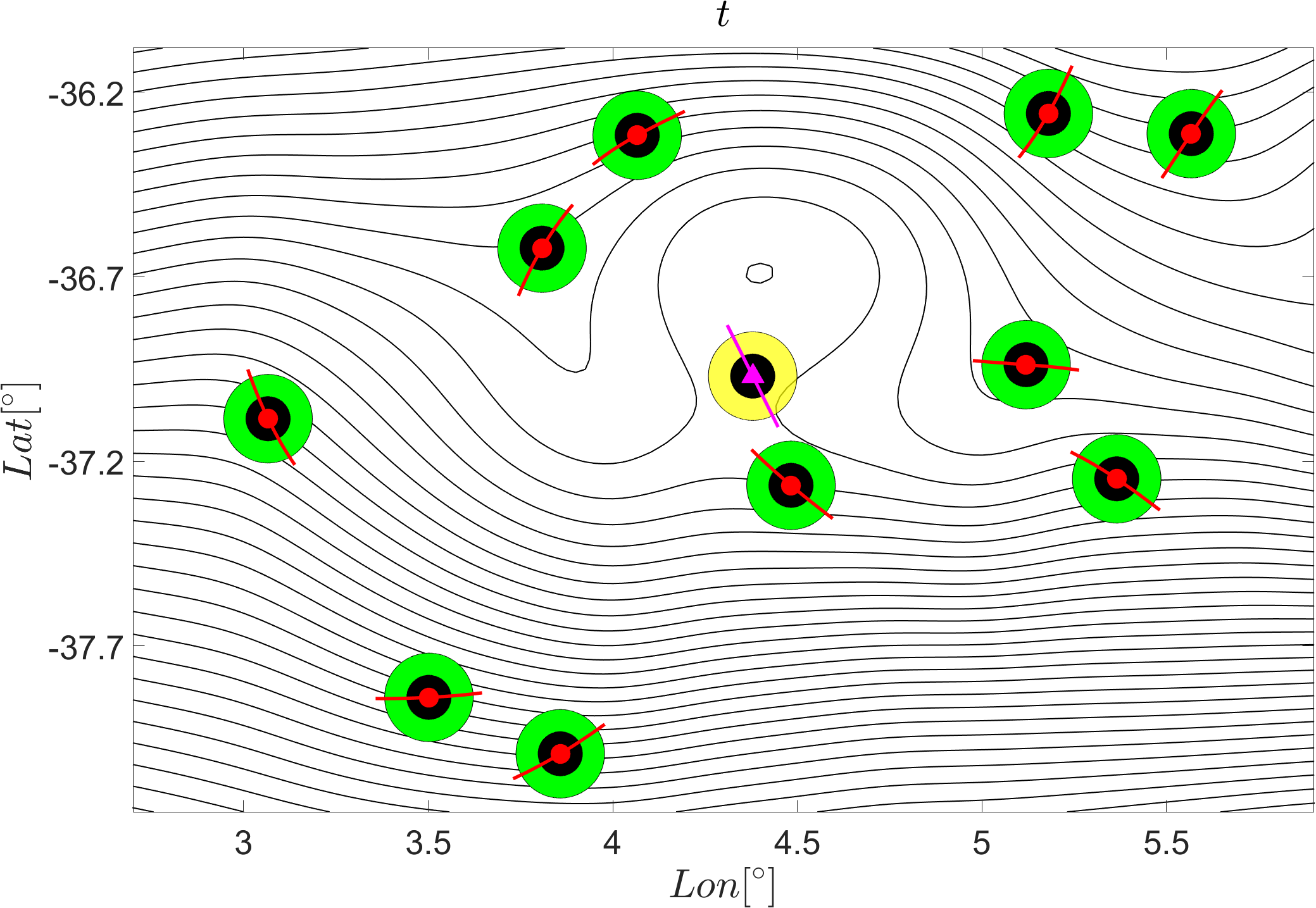}\label{fig:HypOECSsMotMFrt0}}\hfill{}
	\subfloat[]{\includegraphics[width=0.47\textwidth]{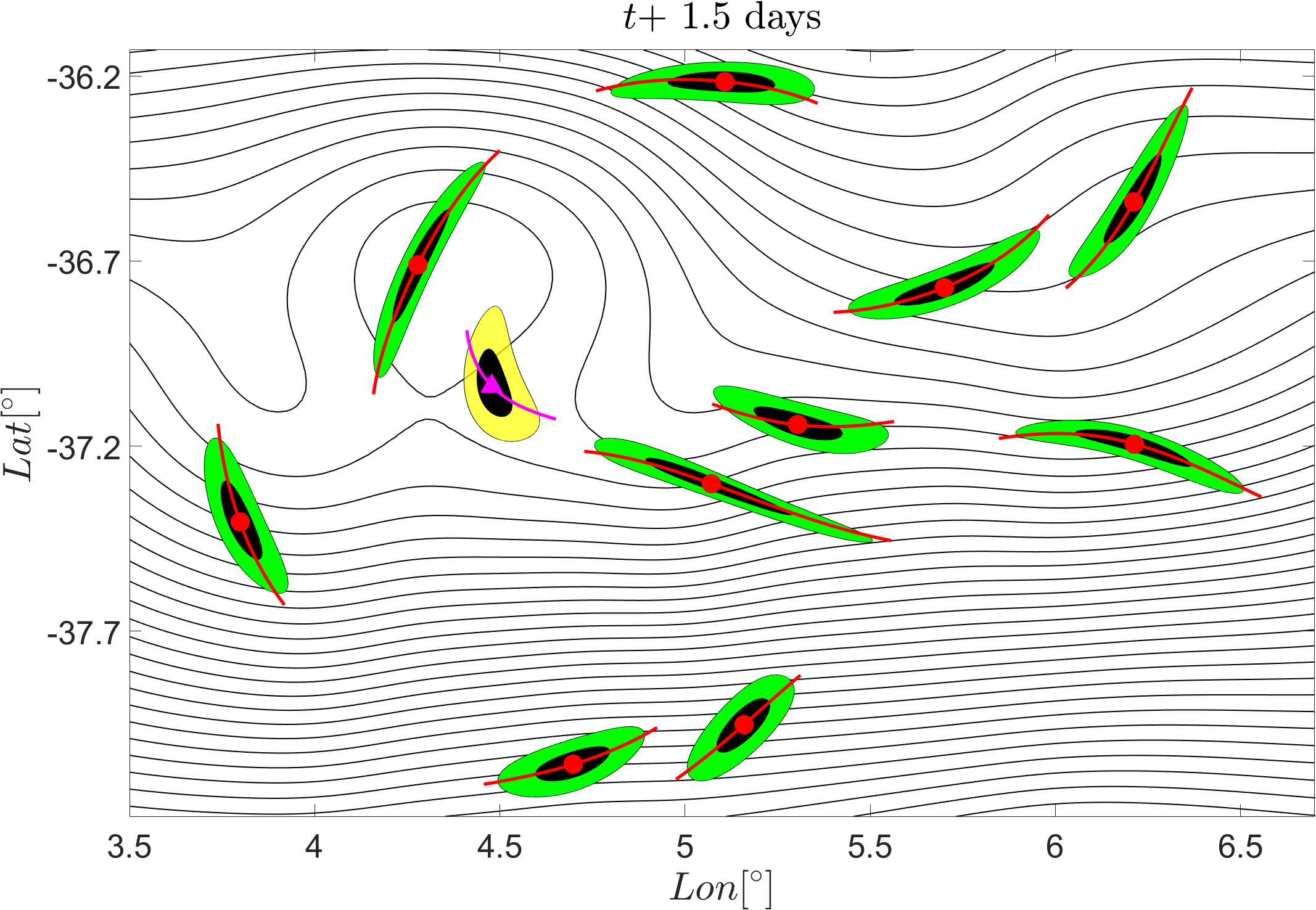}\label{fig:HypOECSsMotMFrtf}}\\
	\caption{(a) Attracting hyperbolic OECSs (red) with their cores (red dots)
		and classic saddle-type stagnation points (magenta triangles) with
		their corresponding unstable directions (magenta), overlaid on streamlines.
		(b) Advected image of the hyperbolic OECSs and classic saddle-type
		stagnation points after 1.5 days, and their effect on nearby particles.
		(c-d) Same phenomenon shown in (a-b), seen from an observer moving
		with constant longitudinal velocity (-0.6 degree/day) relative
		to the one used in (a-b).}	
	\label{fig:HypOECSMotivation}
\end{figure}

As an illustration of the results we derive here, Fig. \ref{fig:HypOECSMotivation}
compares hyperbolic OECSs and the classic frame-dependent saddle-type
stagnation points computed on an ocean velocity data (see section
\ref{NumScheme} for more detail). The effect of those structures
on nearby particles and the instantaneous streamlines are also shown.
In particular, Figs. \ref{fig:HypOECSsMotOrFrt0}-\ref{fig:HypOECSsMotOrFrtf}
show the attracting hyperbolic OECSs (red) with their cores (red dots)
and the saddle-type stagnation points (magenta triangles) with their
unstable directions (magenta) at the initial time, and after short-term
advection. In contrast, Figs. \ref{fig:HypOECSsMotMFrt0}-\ref{fig:HypOECSsMotMFrtf}
show the same features observed from a frame translating in the longitudinal
direction with constant speed (-0.6 degree/day), relative to the one
used in Figs. \ref{fig:HypOECSsMotOrFrt0}-\ref{fig:HypOECSsMotOrFrtf}. 
This example highlights two important facts. First, frame-dependent
diagnostics, such as instantaneous stagnation points, are unsuitable
for the self-consistent identification of coherent structures: the
three saddle-type stagnation points detected by one observer (Figs.
\ref{fig:HypOECSsMotOrFrt0}-\ref{fig:HypOECSsMotOrFrtf}) disappear
when the same phenomenon is analyzed by another observer in a moving
frame (Figs. \ref{fig:HypOECSsMotMFrt0}-\ref{fig:HypOECSsMotMFrtf}).
In this second frame, a single saddle-type stagnation point emerges
at an unrelated location and induces no notable material stretching. 

Second, and more important, hyperbolic OECSs capture actual saddle-type
material behavior in several locations with no stagnation points in
their vicinities. Similar conclusions hold for elliptic and parabolic
OECSs (cf. section \ref{NumScheme}). Based on this example, we expect
OECSs to be useful in a number of physical situations where reliable
now-casting, or short-term forecasting and control of material transport
is critical. 

We develop a theoretical foundation for our global variational OECS
theory in sections 1-5. Readers interested mainly in a practical implementation
of automated OECSs detection may proceed directly to section 6, which
gives a summary of the results along with corresponding numerical
schemes. In section 7, we perform OECSs detection in a two-dimensional
ocean velocity data obtained from satellite altimetry. We show that
variational OECSs outperform traditional Eulerian diagnostics in locating
the skeletons of short-term material deformation.
\section{Set-up and notation}

Consider the two-dimensional non-autonomous dynamical system 
\begin{equation}
\dot{x}=v(x,t),\label{eq:FlowODE}
\end{equation}
with a twice continuously differentiable velocity field $v(x,t)$
defined over an open flow domain $U\in\mathbb{R}^{2}$, over a time
interval $t\in[a,b].$ We recall the customary velocity gradient decomposition
\[
\nabla v(x,t)=S(x,t)+W(x,t),
\]
with the rate-of-strain tensor $S=\tfrac{1}{2}(\nabla v+\nabla v^{T})$
and the spin tensor $W=\tfrac{1}{2}(\nabla v-\nabla v^{T})$. By our
assumptions, $S(x,t)$ and $W(x,t)$are
continuously differentiable in $x$ and $t.$ Under an observer change
\eqref{eq:observer change}, the new rate-of-strain tensor and the
new spin tensor are obtained in the form 
\begin{align}
	\tilde{S}(\tilde{x},t) & =Q^{T}(t)S(x,t)Q(t),\label{eq:strain object}\\
	\tilde{W}(\tilde{x},t) & =Q^{T}(t)W(x,t)Q(t)-Q^{T}(t)\dot{Q}(t),\nonumber 
\end{align}
as shown is classic texts on continuum mechanics (see, e.g., \cite{TruesdellNoll2004}).
Therefore, the rate-of-strain tensor is objective, as it transforms
as a linear operator, whereas the spin tensor is not objective.

The eigenvalues $s_{i}(x,t)$ and eigenvectors $e_{i}(x,t)$ of $S(x,t)$
are defined, indexed and oriented here through the relationship
\begin{equation}
Se_{i}=s_{i}e_{i},\qquad\left|e_{i}\right|=1,\ \ i=1,2;\ \ s_{1}\leq s_{2},\quad e_{2}=Re_{1},\qquad R:=\left(\begin{array}{cc}
0 & -1\\
1 & 0
\end{array}\right).\label{eq:RelativeOrientation}
\end{equation}
We also recall that the rate of length change for an infinitesimal
material element vector $\ell$ based at $x$ is 
\begin{equation}
\dfrac{1}{2}\dfrac{d}{dt}\left|\ell\right|^{2}=\langle\ell,S(x,t)\ell\rangle.\label{eq:rate of length change}
\end{equation}
A further key relationship between the flow map $F_{t_{0}}^{t}:x_{0}\mapsto x(t;t_{0},x_{0})$
of \eqref{eq:FlowODE} and $S(x,t)$ is obtained by considering the
right Cauchy--Green strain tensor
\begin{equation}
C_{t_{0}}^{t}=\left[\nabla F_{t_{0}}^{t}\right]^{T}\nabla F_{t_{0}}^{t},\label{eq:CGdef}
\end{equation}
whose temporal Taylor expansion around the initial time can be computed
as
\begin{equation}
C_{t_{0}}^{t}(x_{0})=I+2S(x_{0},t_{0})(t-t_{0})+\mathcal{\mathcal{O}}\text{}\left(\left|t-t_{0}\right|^{2}\right).\label{eq:CG_Taylor}
\end{equation}
In other words, for small enough times, the leading order Lagrangian
deformation is governed by the Eulerian rate-of-strain tensor. This
observation enables us to consider Eulerian coherent structures as
short-time limits of Lagrangian coherent structures.

\section{Objective deformation rates}

At an arbitrary time $t$, consider a smooth curve $\gamma\subset U$,
parametrized in the form $x(s)$ by its arclength $s\in[0,\sigma]$.
Then the unit vectors $x'(s)$ and $n(s)=R\tfrac{x'(s)}{\left|x'(s)\right|}$,
with the rotation matrix $R$ appearing in \eqref{eq:RelativeOrientation},
define a local tangent and a local normal to $\gamma$ at the point
$x(s)$. The following definition fixes the notions of instantaneous
material shear rate and material stretching rate along $\gamma.$
While these quantities are intuitively clear, we also give their detailed
derivation in Appendix A as first-order terms in the temporal Taylor
expansion of the analogous finite-time Lagrangian shear and stretching
measures (cf. \cite{LCSHallerAnnRev2015}).
\begin{defn}
	Material deformation rates at time $t$ along a material curve $\gamma$
	with arclength parametrization $x(s)$: 
	\begin{enumerate}
		\item \emph{Material stretching rate:}
		\[
		\dot{q}(x,x',t)=\dfrac{\langle x',S(x,t)x'\rangle}{\langle x',x'\rangle}
		\]
		
		\item \emph{Material shear rate:} 
		\begin{equation}
		\dot{p}(x,x',t)=\dfrac{\left\langle x',[S(x,t)R-RS(x,t)]x'\right\rangle }{\left\langle x',x'\right\rangle }.\label{eq:pdotdef}
		\end{equation}
		
	\end{enumerate}
\end{defn}
Physically, $\dot{q}$ gives the instantaneous tangential stretching
rate along $\gamma$, while $\dot{p}$ represents twice the rotation
rate due to shear of an initially normal perturbation to $\gamma$.
The objectivity of these Eulerian deformation rate measures follows
from \eqref{eq:strain object} together with the commutation between
an arbitrary planar rotation $Q(t)$ and the ninety-degree rotation
$R$.

\section{Variational principles for OECS\lowercase{s}}
Using the Eulerian deformation rates introduced in Definition 1, we
now define the averaged stretch- and shear-rate functionals over an
arbitrary curve $\gamma$ at time $t$: 
\begin{defn}
	At time $t$ along the material curve $\gamma$: 
	\begin{enumerate}
		\item The \emph{averaged material stretch-rate} is 
		\begin{equation}
		\dot{Q}_{t}(\gamma)=\dfrac{1}{\sigma}\int_{\gamma}\dot{q}(x(s),x^{\prime}(s),t)\,ds.\label{eq:StretchEAV}
		\end{equation}
		
		\item The \emph{averaged material shear-rate} is 
		\begin{equation}
		\dot{P}_{t}(\gamma)=\dfrac{1}{\sigma}\int_{\gamma}\dot{p}(x(s),x^{\prime}(s),t)\,ds.\label{eq:ShearEAV}
		\end{equation}
		
	\end{enumerate}
\end{defn}
By smooth dependence of $v(x,t)$ on $x$ and $t$, one expects to
see $\mathcal{O}(\epsilon)$ variability in the average material stretch-
or shear-rates across an $\mathcal{O}(\epsilon)$ strip of material
lines surrounding $\gamma.$ Exceptional choices of $\gamma$, however,
defy this trend, displaying only an $\mathcal{O}(\epsilon^{2})$ variability
of $\dot{Q}_{t}(\gamma)$ or $\dot{P}_{t}(\gamma)$ in $\mathcal{O}(\epsilon)$
strips around $\gamma.$ Such $\gamma$ curves will act as centerpieces
of coherence in the material stretch-rate or material shear-rate fields.

This lack of leading order variability in the averaged stretch-rate
or shear-rate along $\gamma$ is equivalent to the vanishing of the
first variation of the functional $\dot{Q}_{t}$ or $\dot{P}_{t}$
on $\gamma$. The former case occurs on instantaneous vortex boundaries
with no short-term unevenness in their tangential deformation (no
short-term filamentation). The latter case arises for curves that
are cores of short-term shear-type deformation (instantaneous jets)
or short-term hyperbolic stretching (instantaneous hyperbolic structures).
We formalize these definitions below, using the time derivatives of
the corresponding global variational LCS definitions reviewed in \cite{LCSHallerAnnRev2015}.
\begin{defn}
	At time $t$ 
	\begin{enumerate}
		\item \textit{A closed curve }$\gamma$ \emph{is an elliptic OECS if it
			is a stationary curve of the averaged stretch-rate functional $\dot{Q}$,
			i.e.,} 
		\begin{equation}
		\delta\dot{Q}_{t}(\gamma)=0.\label{eq:StretchEAV-1}
		\end{equation}
		
		\item \emph{A curve $\gamma$ is a shearless} OECS \emph{if it is a stationary
			curve of the averaged shear-rate functional $\dot{P}_{t}$, i.e.,}
		\begin{equation}
		\delta\dot{P}_{t}(\gamma)=0.\label{eq:ShearEAV-1}
		\end{equation}
		
	\end{enumerate}
\end{defn}
The variational problems outlined in \eqref{eq:StretchEAV-1} and
\eqref{eq:ShearEAV-1} are equivalent to the weak Euler-Lagrange equations
\begin{equation}
\delta\dot{Q}_{t}(\gamma)=\dfrac{1}{\sigma}[\langle\partial_{x'}\dot{q},h\rangle]_{0}^{\sigma}+\dfrac{1}{\sigma}\int_{0}^{\sigma}\left[\partial_{x}\dot{q}-\dfrac{d}{ds}\partial_{x'}\dot{q}\right]h\,ds=0,\label{eq:AvStrRatFunEulEq}
\end{equation}
\begin{equation}
\delta\dot{P}_{t}(\gamma)=\dfrac{1}{\sigma}[\langle\partial_{x'}\dot{p},h\rangle]_{0}^{\sigma}+\dfrac{1}{\sigma}\int_{0}^{\sigma}\left[\partial_{x}\dot{p}-\dfrac{d}{ds}\partial_{x'}\dot{p}\right]h\,ds=0,\label{eq:AvShRatFunEulEq}
\end{equation}
with $h(s)$ denoting small perturbations to the curve $\gamma$.
We discuss below the solutions of equations (\ref{eq:AvStrRatFunEulEq}-\ref{eq:AvShRatFunEulEq}),
building on similar ideas developed for LCSs (\cite{BlackHoleHaller2013,ShearlessBarrierFarazmand2014}).
This approach will lead to explicit solutions of\eqref{eq:StretchEAV-1}-\eqref{eq:ShearEAV-1},
and enable a further partitioning of shearless OECSs into hyperbolic
and parabolic OECSs.

\section{Elliptic OECS}

By Definition 3, an elliptic OECS is a closed curve, and hence its
small perturbation $h(s)$ used in calculating the first variation
of $\dot{Q}_{t}$ along $\gamma$ is also periodic with the same period
$\sigma$. As a consequence, the first bracketed term in \eqref{eq:AvStrRatFunEulEq}
vanishes. Then, by the fundamental lemma of the calculus of variations,
equation \eqref{eq:StretchEAV-1} becomes equivalent to the classic
Euler--Lagrange equations
\begin{equation}
\partial_{x}\dot{q}-\dfrac{d}{ds}\partial_{x'}\dot{q}=0.\label{eq:EL-stretch}
\end{equation}

The detailed form of this equation is similar to the Euler--Lagrange
equation derived by \cite{BlackHoleHaller2013} for elliptic LCSs
(stationary curves of the averaged, finite-time tangential stretch
along the material curve $\gamma$). The only difference in the present
context is that the expression for $\dot{q}$ has no square root,
and contains the rate-of-strain tensor $S(x,t)$ instead of the Cauchy--Green
strain tensor $C_{t_{0}}^{t}(x_{0})$ defined in \eqref{eq:CGdef}.
Therefore, following the procedure similar to the one adopted in \cite{BlackHoleHaller2013},
we obtain the following full characterization of elliptic OECSs.
\begin{thm}
	Elliptic OECSs at time $t$ coincide with limit cycles of the differential
	equation family 
	\begin{equation}
	\frac{dx}{ds}=\chi_{\mu}^{\pm}(x,t),\qquad\chi_{\mu}^{\pm}=\sqrt{\dfrac{s_{2}-\mu}{s_{2}-s_{1}}}e_{1}\ \pm\ \sqrt{\dfrac{\mu-s_{1}}{s_{2}-s_{1}}}e_{2},\label{eq:chimufield}
	\end{equation}
	defined for any parameter value $\mu\in\mathbb{R}$ on the set 
	\[
	U_{\mu}:=\{x\in U:\ s_{1}\neq s_{2},\ s_{1}\le\mu\le s_{2}\}.
	\]
	Along any elliptic OECS obtained for a given value of $\mu$, the
	material stretching rate at time $t$ is pointwise constant and equal
	to $\dot{q}\equiv\mu.$ 
\end{thm}
The right-hand side of the differential equation \eqref{eq:chimufield}
is only locally a vector field in $U_{\mu}$, because the rate-of-strain
eigenvector fields $e_{i}$ are generally not globally orientable
in $U$. Local orientability of $\chi_{\mu}^{\pm}(x,t)$ in $U_{\mu}$,
however, is always possible. This makes the trajectories, and specifically
the limit cycles, of the differential equation \eqref{eq:chimufield}
well defined. 

Along limit cycles obtained for $\mu=0$, we have $\dot{q}\equiv0.$
Such elliptic OECSs, therefore, are perfectly coherent in the Eulerian
sense, exhibiting uniformly zero pointwise material stretching rate.
The direction field $\chi_{\mu}^{\pm}$ ceases to be well-defined
at locations of repeated rate-of-strain eigenvalues ($s_{1}=s_{2}$).
Following the terminology used in \cite{BlackHoleHaller2013}, we
refer to such locations as \emph{singularities} of the tensor field
$S(\,\cdot\,,t).$ An appropriate extension of Poincaré's classic
index theory implies that there are at least two singularities of
$S$ in the interior of any limit cycle of $\chi_{\mu}^{\pm}$. This
in turn leads to an automated detection algorithm that applies equally
to elliptic OECSs. We refer to \cite{AutoMdetectFlorian} for a detailed
discussion of this algorithm for elliptic LCSs.

\begin{figure}[h]
	\subfloat[\label{fig:chi0Stretch}]{\includegraphics[width=0.45\textwidth]{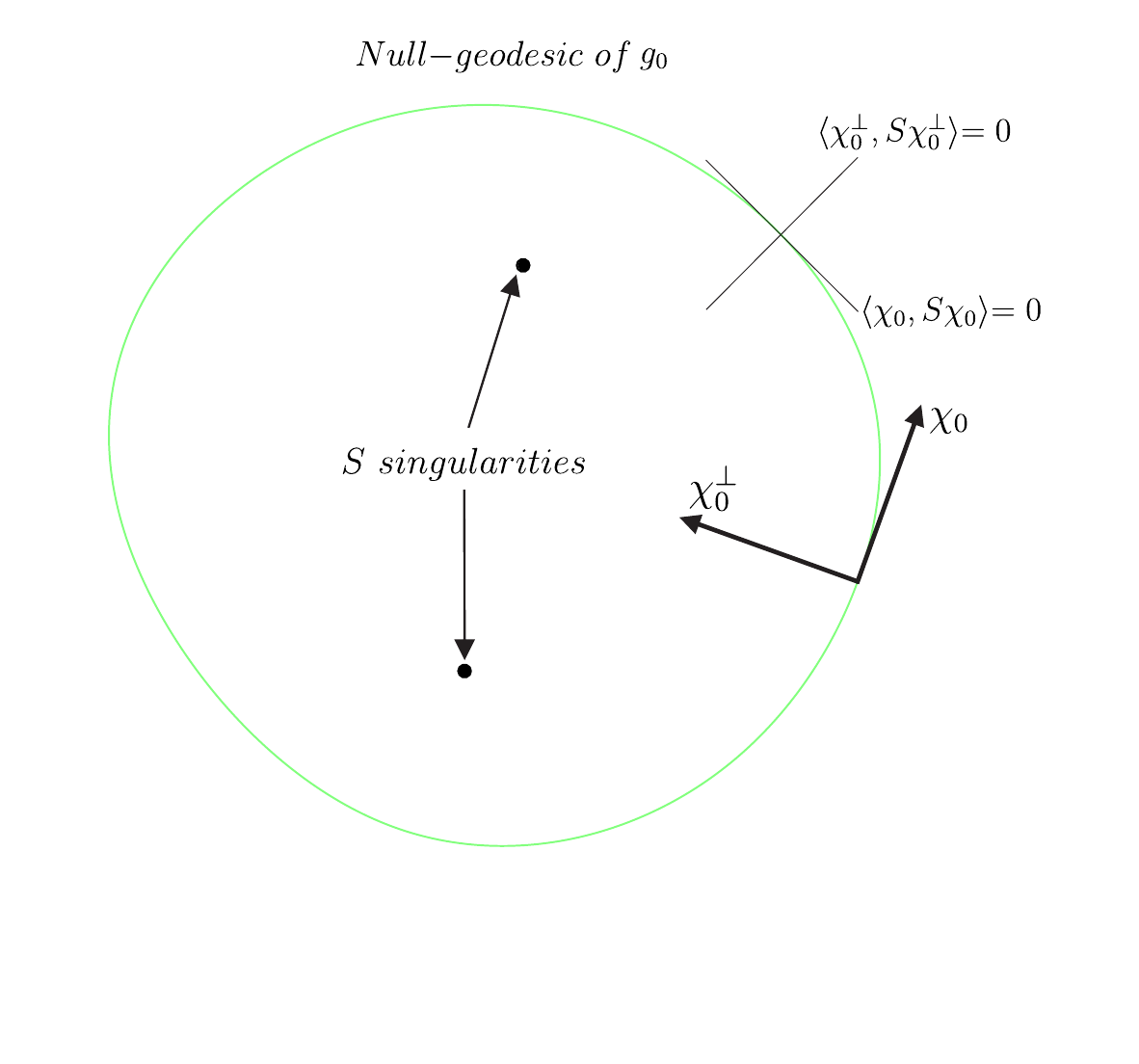}}\hfill{}
	\subfloat[\label{fig:BeltofellipticECSs}]{\includegraphics[width=0.5\textwidth]{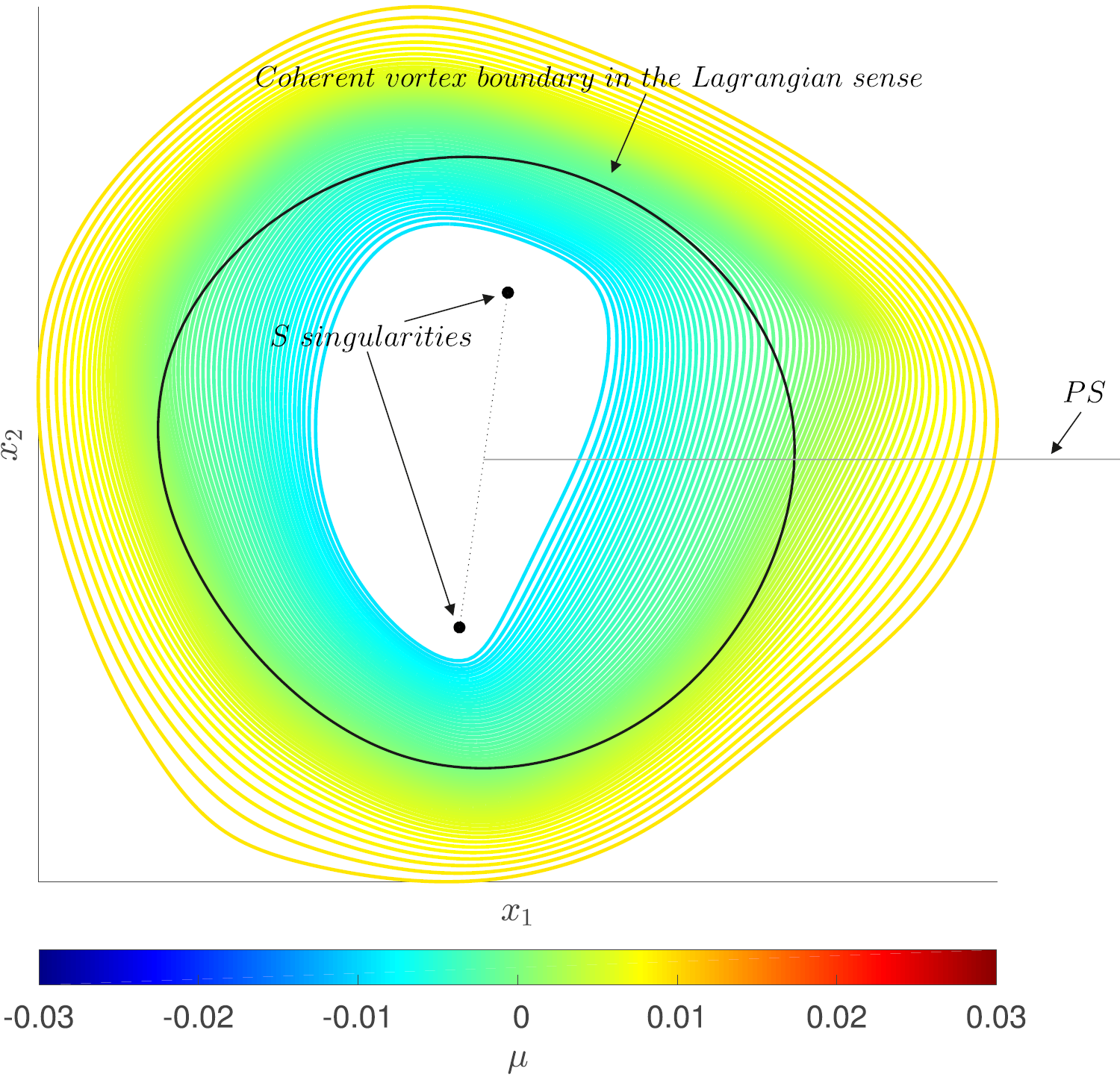}}\\
	\caption{(a) Perfectly coherent elliptic OECS at time $t$, obtained as a limit cycle $\gamma$ of the direction field $\chi_{0}^{\pm}(x,t)$. The tangential stretching rate (and, in incompressible flows, the normal repulsion rate) along $\gamma$ is pointwise zero. Furthermore, the interior of $\gamma$ always contains at least two singularities of the tensor field $S$. (b) Nested family of elliptic OECSs for different values of $\mu$ (in color). A nearby Lagrangian vortex boundary is shown in black. Black dots represent singularities of the rate-of-strain tensor used for an automatic placement of the Poincaré Section (PS) in the detection of limit cycles for the direction field family $\chi_{\mu}^{\pm}(x,t)$.}	
	\label{fig:eta0_ISGs}
\end{figure}

By their hyperbolicity (i.e., strict attraction or repulsion of nearby
trajectories of the vector field $\chi_{\mu}^{\pm}$), elliptic OECSs
are robust with respect to moderate errors and uncertainties in the
underlying velocity field $v(x,t)$. Figure \ref{fig:chi0Stretch}
summarizes the main properties of perfectly coherent elliptic OECSs.
Detailed discussions on the numerical detection of limit cycles of
directions fields can be found in \cite{BlackHoleHaller2013} and\textbf{
}\cite{AutoMdetectFlorian} for LCS. These can be directly applied
here after replacing $C_{t_{0}}^{t}(x_{0})$ with $S(x,t)$. 

Under changes in the parameter $\mu$, limit cycles of $\chi_{\mu}^{\pm}$
arise in continuous, non-intersecting families (see Appendix B). An
\emph{Eulerian vortex boundary} can then be defined as the outermost
member of such an elliptic OECS family. Given the objectivity of each
member of such a limit cycle family, Eulerian vortex boundaries defined
in this fashion are also objective. 

Figure \ref{fig:BeltofellipticECSs} shows such a family and its outermost
member in a flow example analyzed in more detail in section 10. A
nested OECSs family often signals a nearby coherent Lagrangian vortex
boundary, as is the case in Fig. \ref{fig:BeltofellipticECSs}. Just
as in the case of elliptic LCSs \cite{BlackHoleHaller2013}, solutions
of the variational principle \eqref{eq:StretchEAV-1} can also be
viewed as closed null-geodesics of a Lorentzian metric family 
\begin{equation}
g_{\mu}(u,u)=\langle u,S_{\mu}u\rangle,\ \ \ \ S_{\mu}(x,t)=[S(x,t)-\mu I],\label{eq:LorMetrStetch}
\end{equation}
which has metric signature (-,+) \cite{Beem1996}. The tensor family
$S_{\mu}$ denotes a generalized rate-of-strain tensor and the parameter
$\mu$ denotes the instantaneous tangential stretching rate on $\gamma$.
In this context, locations of repeated eigenvalues of $S(x,t)$ are
singularities for the metric $g_{\mu}$ which becomes degenerate at
these points (cf. Fig. \ref{fig:eta0_ISGs}).

\section{Shearless OECS}

By Definition 3, a shearless OECS at time $t$ is a stationary curve
of the averaged shear-rate functional $\dot{P}_{t}$, and hence satisfies
the weak Euler--Lagrange equation \eqref{eq:AvShRatFunEulEq}. We
can pass to the strong form of these Euler--Lagrange equation if the
first bracketed term in \eqref{eq:AvShRatFunEulEq} vanishes for the
class of admissible perturbations $h(s)$ applied to the parametrization
$x(s)$ of the stationary curve $\gamma$ we seek. 

Following the procedure developed in \cite{ShearlessBarrierFarazmand2014}\textbf{
}for elliptic LCS\textbf{, }we find the following class of admissible
perturbations for which the boundary term $\langle\partial_{x'}\dot{p},h\rangle$
in \eqref{eq:AvShRatFunEulEq} vanishes:
\begin{description}
	\item [{\emph{(BC1)}}] \emph{Variable-endpoint boundary conditions:} The
	endpoints of $\gamma$ coincide with singularities of the rate-of-strain
	tensor field $S(x,t)$, i.e., we have $s_{1}(x(0),t)=s_{2}(x(0),t)$
	and $s_{1}(x(\sigma),t)=s_{2}(x(\sigma),t)$. In this case, the perturbations
	$h(s)$ are arbitrary, including arbitrary perturbations to the endpoints
	of $\gamma$. 
	\item [{\emph{(BC2)}}] \emph{Fixed-endpoint boundary conditions}: The perturbation
	$h(s)$ is arbitrary for all $s\in(0,\sigma),$ but must leave the
	endpoints of $\gamma$ fixed: $h(0)=h(\sigma)=0$.
\end{description}
Under (BC1) or (BC2), the bracketed term in \eqref{eq:AvShRatFunEulEq}
vanishes. By the fundamental lemma of the calculus of variations,
a solution curve $\gamma$ of \eqref{eq:AvShRatFunEulEq} must then
satisfy the Euler--Lagrange equations 

\begin{equation}
\partial_{x}\dot{p}-\dfrac{d}{ds}\partial_{x'}\dot{p}=0.\label{eq:EulLagrEqSh}
\end{equation}
Invoking the results of \cite{ShearlessBarrierFarazmand2014}\textbf{
}for shearless LCSs, we obtain that solution curves of \eqref{eq:EulLagrEqSh}
with exactly vanishing shear rates are piecewise tangent to one of
the eigenvector fields of $S(x,t)$. 
\begin{thm}
	Shearless OECSs at time $t$ coincide with continuous trajectories
	of the differential equation family 
	\begin{equation}
	\frac{dx}{ds}=e_{i}(x,t),\qquad S(x,t)e_{i}(x,t)=s_{i}(x,t)e_{i}(x,t),\qquad i=1,2.\label{eq:e_i field}
	\end{equation}
	Along any such shearless OECS, the pointwise material shear rate
	at time $t$ is zero.
\end{thm}
We refer to the trajectories of \eqref{eq:e_i field} as \emph{$e_{1}$-lines
}and \emph{$e_{2}$-lines}, respectively. We will use the boundary
condition classes (BC1)-(BC2) to further distinguish parabolic and
hyperbolic OECSs within the shearless OECSs satisfying \eqref{eq:e_i field}.

From an argument closely following \cite{ShearlessBarrierFarazmand2014},
we obtain that the solutions of \eqref{eq:EulLagrEqSh} can also be
viewed as null-geodesics of the Lorentzian metric family

\[
m_{\nu}(u,u)=\langle u,[2S(x,t)R-\nu I]u\rangle,\qquad\nu\in\mathbb{R},
\]
which again has metric signature (-,+) \cite{Beem1996}, and admits
singularities at locations of repeated eigenvalues of $S(x,t)$. In
particular, shearless OECSs are null geodesics of $m_{\nu}(u,u)$
for $\nu=0.$

\subsection{Parabolic OECS }

\emph{Parabolic OECSs }are trajectories of \eqref{eq:e_i field} that
satisfy the free-endpoint boundary conditions (BC1) and are as close
as possible to being neutrally stable, as detailed below. By the nature
of (BC1), such OECSs are stationary curves of the shear-rate functional
$\dot{P}_{t}(\gamma)$ under the broadest possible set of perturbations.
This makes parabolic OECSs the most observable class of shearless
OECSs, creating \textit{\emph{short-term pathways}} that mimic the
role of the Lagrangian jet cores identified in\textbf{ }\cite{ShearlessBarrierFarazmand2014}.

More specifically, parabolic OECSs are heteroclinic chains of $e_{1}-$lines
and $e_{2}$-lines connecting singularities of $S(x,t)$. For observability
and uniqueness, we only consider alternating chains of $e_{1}$- and
$e_{2}$-line connections that are locally unique and structurally
stable. As shown in \cite{Delmarcelle1994}, the only structurally
stable tensorline singularities are trisectors and wedges, shown in
Fig. 3. 

\begin{figure}[h]
	\centering
	\includegraphics[width=0.5\textwidth]{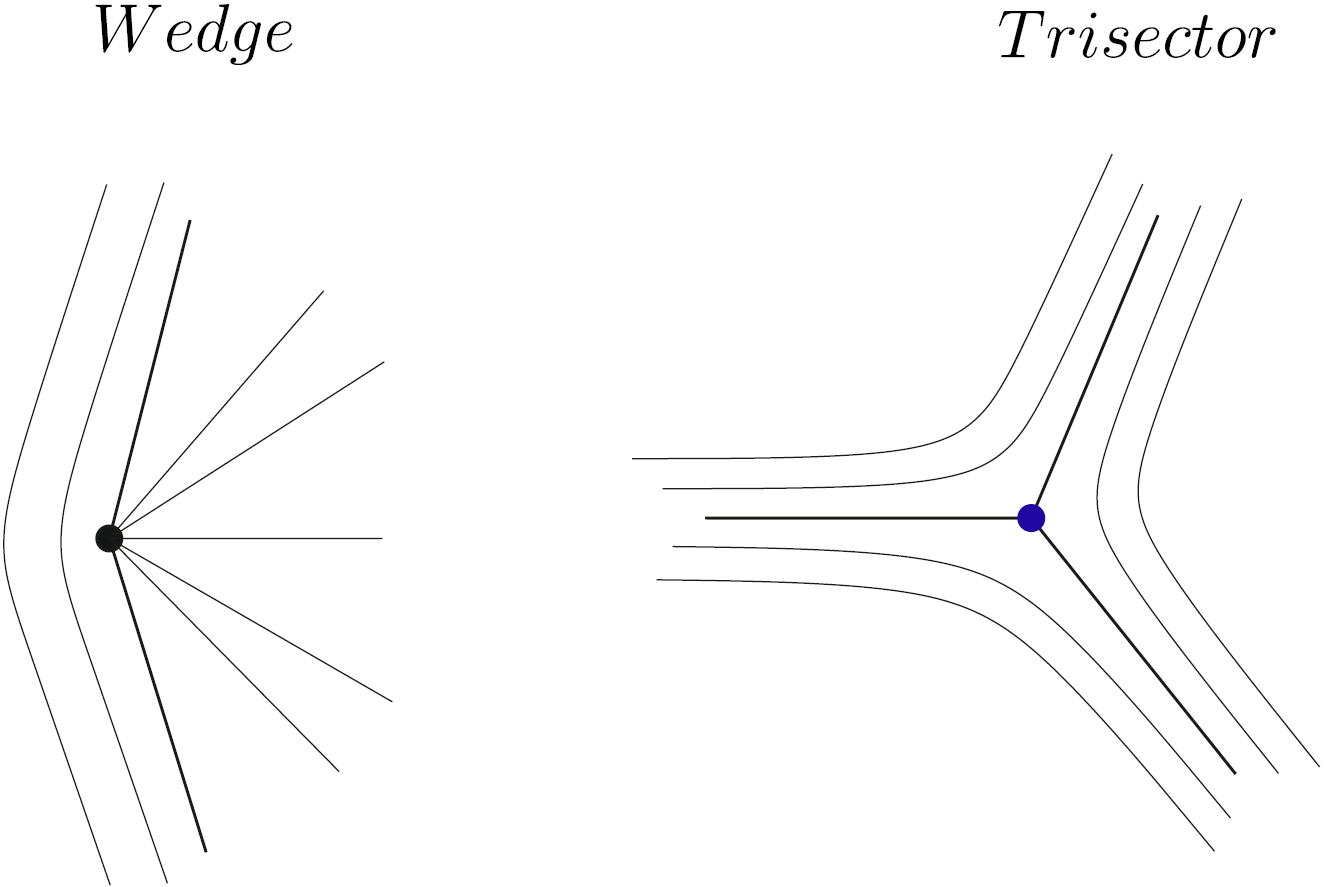}
	\caption{Local topology of a tensorline field around structurally stable singularities:
		Wedge and Trisector.}
	\label{fig:SingTypes}
\end{figure}

Furthermore, as shown in \cite{ShearlessBarrierFarazmand2014}, a
structurally stable and unique connection between two such singularities
must necessarily be a wedge-trisector connection, as illustrated in
Fig. 4.
\begin{figure}[h]
\includegraphics[width=0.9\textwidth]{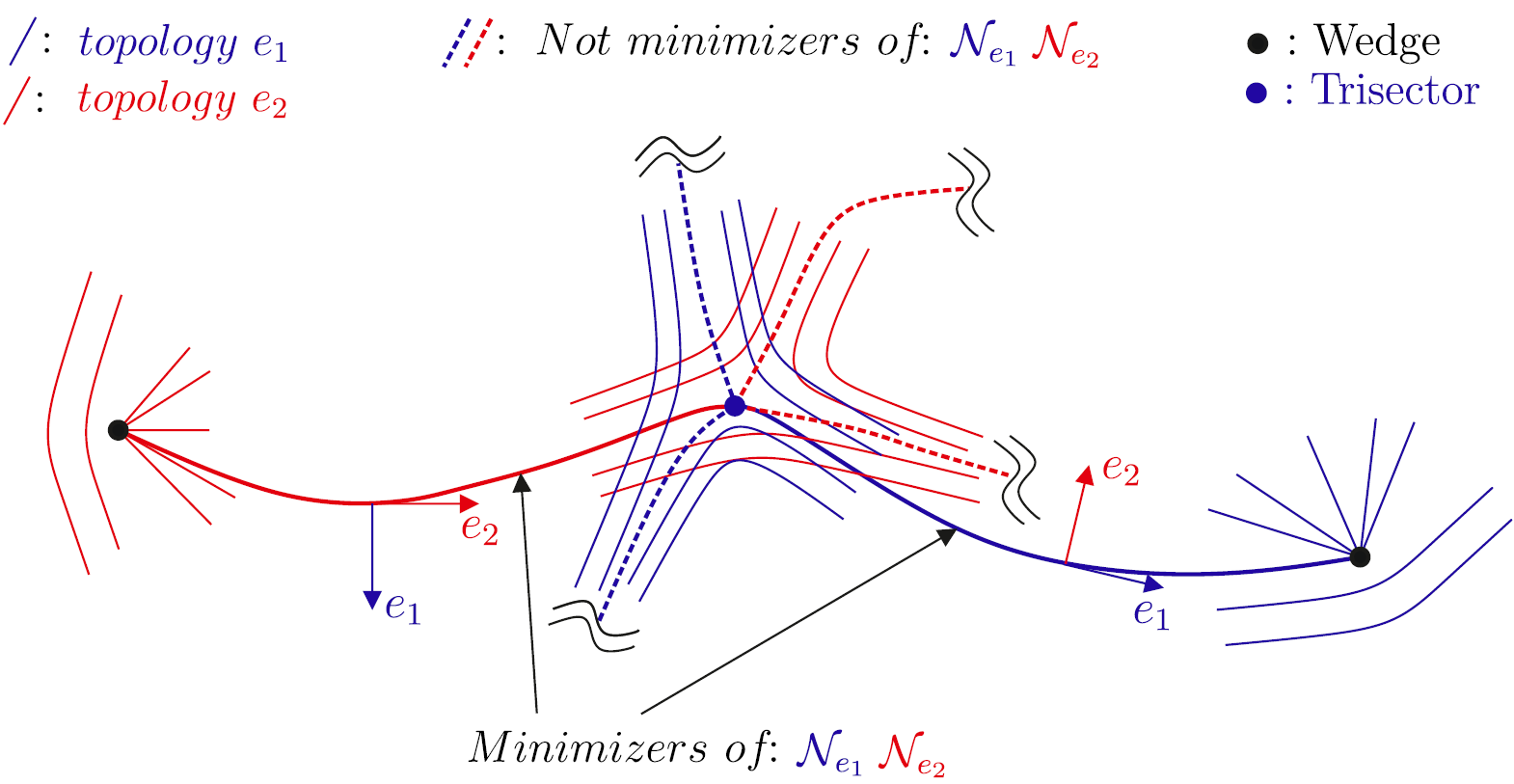}
	\caption{Parabolic OECSs are alternating chains of $e_{i}$-lines that connect
		structurally stable singularities and are weak minimizers of the corresponding
		neutrality functions $\mathcal{N}_{e_{i}}$.}
	\label{fig:ParabolicECSs}
\end{figure}
By definition, each segment of an alternating chain of $e_{1}$- and
$e_{2}$-lines repels or attracts trajectories over short enough time
intervals. To ensure that neither attraction nor repulsion prevails
for the whole trajectory chain, we require the chain to be as close
as possible to being neutrally stable. To this end, we introduce the
pointwise \emph{neutrality functions} 
\begin{equation}
\mathcal{N}_{e_{1}}(x,t)=s_{2}^{2}(x,t),\ \ \mathcal{N}_{e_{2}}(x)=s_{1}^{2}(x,t),\label{eq:Neutralitydef}
\end{equation}
with the function $\mathcal{N}_{e_{i}}(x,t)$ measuring how close
the squared rate of attraction or repulsion, $s_{j}^{2}(x,t),$ along
an $e_{i}$ line is to zero at time $t$. In case of incompressible
flows, we have $\mathcal{N}_{e_{1}}(x,t)=\mathcal{N}_{e_{2}}(x,t)$,
given that $s_{1}(x,t)+s_{2}(x,t)=0$. 

All this follows closely the variational approach developed in \cite{ShearlessBarrierFarazmand2014}
for parabolic LCSs, but with $S(x,t)$ substituted for $C_{t_{0}}^{t}(x_{0})$.
As a next step, we introduce the convexity sets 
\begin{equation}
\mathcal{C}_{e_{i}}(t)=\{x\in U:\langle e_{j}(x,t),\partial_{r}^{2}\mathcal{N}_{e_{i}}(x,t)e_{j}(x,t)\rangle>0,\ \ i\neq j\},\ \ i=1,2.\label{eq:Convexityset}
\end{equation}
Each such set $\mathcal{C}_{e_{i}}(t)$ is simply the set of points
at which the corresponding neutrality $\mathcal{N}_{e_{i}}(x,t)$
is a convex function of $x$ at time $t$. We say that a compact $e_{i}$-line
segment $\gamma$ is a \emph{weak minimizer} of $\mathcal{N}_{e_{i}}$
at time $t$ if both $\gamma$ and the nearest trench of the function
$\mathcal{N}_{e_{i}}(\,\cdot\,,t)$ lie in the same connected component
of $\mathcal{C}_{e_{i}}(t)$. More precisely, if the arclength parametrization
of $\gamma$ is $x(s)$, and the unit normal along $\gamma$ is given
by $n(s)$, then we require 

\begin{equation}
x(s)+\epsilon n(s)\in\mathcal{C}_{e_{i}}(t),\ \ s\in[0,\sigma],\quad\epsilon\in[0,\epsilon_{0}(s,t)],\label{eq:WeakminimofNeutrality}
\end{equation}
where 
\begin{equation}
\epsilon_{0}(s,t)=\min\{\lvert\epsilon\rvert\in\mathbb{R}^{+}:\partial_{\epsilon}\mathcal{N}_{e_{i}}(x(s)+\epsilon n(s),t)=0,\ \ \partial_{\epsilon}^{2}\mathcal{N}_{e_{i}}(x(s)+\epsilon n(s),t)>0\}.\label{eq:WeakminimofNeutrality2}
\end{equation}

We summarize our formal definition of parabolic OECSs as follows (cf.
Fig. 4). 
\begin{defn}
	A \emph{parabolic OECS} at time $t$ is a shearless OECS composed
	of alternating chains of $e_{1}$- and $e_{2}$-line segments that
	connect wedge and trisector singularities of the the rate-of-strain
	tensor $S(x,t).$ Furthermore, each $e_{i}$ segment in the chain
	is a weak minimizer of the neutrality function $\mathcal{N}_{e_{i}}(x,t).$ 
\end{defn}

\subsection{Hyperbolic OECSs}

\emph{Hyperbolic OECSs }are trajectory segments of \eqref{eq:e_i field}
that satisfy the fixed-endpoint boundary conditions (BC2) and contain
precisely one point of maximal repulsion-rate or maximal attraction-rate.
This point will then play the role of an instantaneous saddle point,
with the OECSs acting as the short-term stable or unstable manifold
for this saddle point. By the nature of (BC2), hyperbolic OECSs are
only stationary curves of the shear-rate functional $\dot{P}_{t}(\gamma)$
under variations that leave their endpoints fixed. This makes individual
hyperbolic OECSs less observable than parabolic OECSs. This is also
the case for hyperbolic LCSs, which are generally responsible for
the intricate filamentation of tracer patterns. Details of these filaments
are generally less observable and robust as those of Lagrangian jet
cores marked by parabolic LCSs (cf. \cite{ShearlessBarrierFarazmand2014}).

Within the family of hyperbolic OECSs, we distinguish\emph{ attracting
	OECSs} as material curves that attract nearby material curves instantaneously.
Similarly, we distinguish \emph{repelling OECSs }as hyperbolic OECSs
that instantaneously repel all nearby material curves. We summarize
these definitions more formally as follows (cf. Fig. 5).
\begin{defn}
	A\emph{ repelling OECS} at time $t$ is an open $e_{1}$-line segment
	that contains a local maximum of the function $s_{2}(\,\cdot\,,t)$,
	but contains no other local extremum point of $s_{2}(\,\cdot\,,t)$.
	An\emph{ attracting OECS} at time $t$ is an open $e_{2}$-line segment
	that contains a local minimum of the function $s_{1}(\,\cdot\,,t)$,
	but contains no other local extremum of $s_{1}(\,\cdot\,,t)$. Finally,
	a \emph{hyperbolic OECS} is a shearless OECS that is either an attracting
	or a repelling OECS. 
\end{defn}
\begin{figure}[h]
	\centering
	\includegraphics[width=0.9\textwidth]{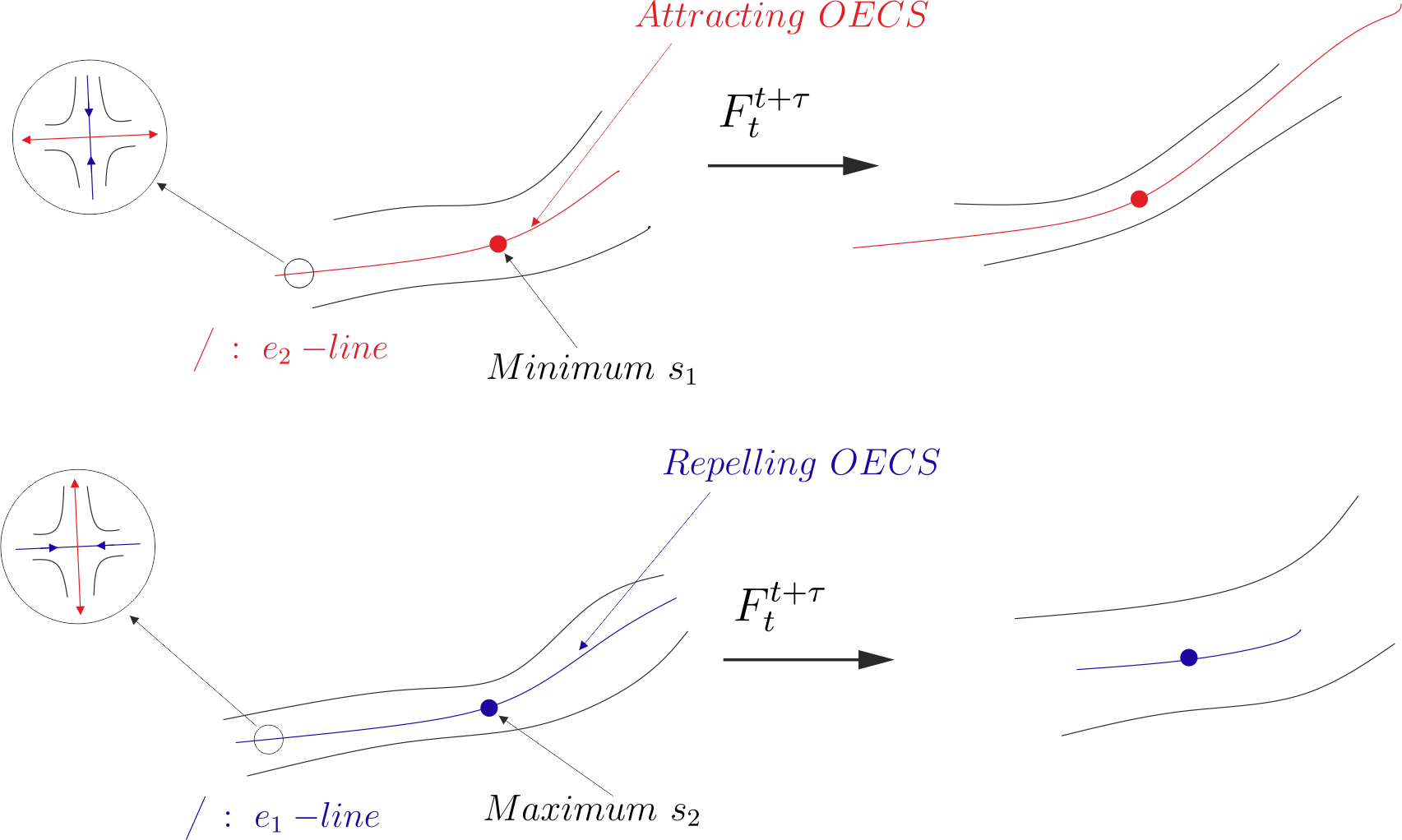}	
	\caption{Attracting (repelling) hyperbolic OECSs as the instantaneous attracting
		(repelling) material lines launched from a minimum of $s_{1}$, (maximum
		of $s_{2}$). In the circular inset: the local tangential stretching
		and normal repulsion.}
	\label{fig:HypECSs}
\end{figure}

\begin{rem}
	As indicated in Fig. 5, the cores of hyperbolic OECSs are defined
	by a local maximum of $s_{2}(\,\cdot\,,t)$ along repelling OECSs
	and by a local minimum of $s_{1}(\,\cdot\,,t)$ along attracting OECSs.
	We call these cores \emph{objective saddle points,} as they represent
	generalizations of classic saddle-type stagnation points from steady
	flows.\textbf{ }Instantaneous stagnation points of an unsteady velocity
	field are not Galilean invariant, i.e., may disappear even under constant-speed
	translations of the coordinate frame (cf. Fig. \ref{fig:HypOECSMotivation}).
	In contrast, the objective saddle points introduced here are objective,
	and hence persist even under general observer changes of the form
	\eqref{eq:observer change}.\textbf{ }The local value of $s_{i}(x,t)$
	on the objective saddle quantifies the strength of that saddle objectively\emph{.}
\end{rem}

\section{Summary of OECS\lowercase{s} and their numerical identification\label{SecAutomScheme} }

In Table \ref{SummaryTableECSs}, we list the OECSs we introduced
in the previous section, together with the ODEs and boundary conditions
they satisfy. 

\begin{table}[h]
	\begin{centering}
		\begin{tabular}{l|l|l}
			\hline 
			$\mathbf{Type\ of\ OECS}$  & $\mathbf{ODE}$  & $\mathbf{Boundary\ condition}$s\tabularnewline
			\hline 
			Attracting  & $x'=e_{2}(x,t)$  & arbitrary; $x(s)$ contains a local minimum of $s_{1}(\,\cdot\,,t)$\tabularnewline
			\hline 
			Repelling & $x'=e_{1}(x,t)$  & arbitrary; $x(s)$ contains a local maximum of $s_{2}(\,\cdot\,,t)$\tabularnewline
			\hline 
			Parabolic  & $x'=e_{i}(x,t)$ with alternating $i$  & $s_{2}(x,t)=s_{1}(x,t)$\tabularnewline
			\hline 
			Elliptic  & $x'=\chi_{\mu}^{\pm}(x,t)$  & Periodic\tabularnewline
			\hline 
		\end{tabular}
		\par\end{centering}
	
	\caption{Summary of the different types of OECSs.}

	\centering{}\label{SummaryTableECSs} 
\end{table}

Next, we summarize the numerical steps in locating OECSs in a planar
unsteady flow. We start with the common steps, then detail the further
steps for different types of OECSs separately. 

\begin{algorithm}[H]
	\protect\caption{Compute the rate-of-strain tensor $S(x,t)$, its invariants and singularities }
	\label{algorithm1} \textbf{Input:} A 2-dimensional velocity field
	$v(x,t)$
	\begin{enumerate}
		\item Compute the rate-of-strain tensor $S(x,t)=\frac{1}{2}\left(\nabla v(x,t)+\left[\nabla v(x,t)\right]^{T}\right)$
		at the current time $t$ on a rectangular grid over the $(x_{1},x_{2}$)
		coordinates. 
		\item Detect the singularities of $S$ as common, transverse zeros of $S_{11}(\,\cdot\,,t)-S_{22}(\,\cdot\,,t)$
		and $S_{12}(\,\cdot\,,t)$, with $S_{ij}$ denoting the entry of $S$
		at row $i$ and column $j$. 
		\item Determine the type of the singularity (trisector or wedge) as described
		in \cite{ShearlessBarrierFarazmand2014}. 
		\item Compute the eigenvalue fields $s_{1}(x,t)<s_{2}(x,t)$  and the associated
		unit eigenvector fields $e_{i}(x,t)$ of $S(x,t)$ for $i=1,2.$
	\end{enumerate}
	\textbf{Output:} $S(x,t)$ as well as $s_{i}(x,t)$ and $e_{i}(x,t)$,
	for $i=1,2$, and the position and type (wedge or trisector) of the
	rate-of-strain singularities $x^{j}(t)$ satisfying $s_{1}(x^{j}(t),t)=s_{2}(x^{j}(t),t),$
	$j=1,\ldots,N.$
\end{algorithm}

\subsection{Elliptic OECSs}

To locate elliptic OECSs automatically as limit cycles of the direction
field \eqref{eq:chimufield}, we rely on a version of Poincaré's index
theory extended to direction fields \cite{DELMARCELLEPhDThesis1994,AutoMdetectFlorian}.
A consequence of this theory is that at least two wedge-type singularities
of $S(x,t)$ must exist inside any such limit cycle \cite{AutoMdetectFlorian}.
For robust limit cycles, we seek nearby wedge pairs surrounded by
an annular region of no singularities. The same procedure arises in
elliptic LCS detection, involving the location of singularities of
the tensor field $C_{t_{0}}^{t_{1}}.$ We refer to \cite{AutoMdetectFlorian}
for details of this numerical algorithm. 
\begin{algorithm}[H]
\caption{Compute Elliptic OECSs}
	\label{algorithm2} \textbf{Input:} $S(x,t)$ as well as $s_{i}(x,t)$
	and $e_{i}(x,t)$, for $i=1,2$, and the position and type (wedge
	or trisector) of the rate-of-strain singularities $x^{j}(t)$ satisfying
	$s_{1}(x^{j}(t),t)=s_{2}(x^{j}(t),t),$ $j=1,\ldots,N.$
	\begin{enumerate}
		\item Locate isolated wedge-type pairs of singularities and place the Poincaré
		sections at their midpoint. 
		\item Compute the vector field $\chi_{\mu}^{\pm}(x(s),t)$ defined in \eqref{eq:chimufield}
		for different values of stretching rate $\mu$, remaining in the range
		$\mu\approx0$. 
		\item Use the Poincaré sections as sets of initial conditions in the computation
		of limit cycles of 
		\[
		x'(s)=\mathrm{sign}\left\langle \chi_{\mu}^{\pm}(x(s)),\tfrac{dx(s-\Delta)}{ds}\right\rangle \chi_{\mu}^{\pm}(x(s)),
		\]
		where the factor multiplying $\chi_{\mu}^{\pm}(x(s),t)$ removes potential
		orientational discontinuities in the direction field $\chi_{\mu}^{\pm}(x(s),t)$
		away from singularities, and $\Delta$ denotes the integration step
		in the independent variable $s$. 
	\end{enumerate}
	\textbf{Output: }Elliptic OECSs. 
\end{algorithm}

\subsection{Hyperbolic OECSs}

\begin{algorithm}[H]
\caption{Compute Hyperbolic OECSs}
	\label{algorithm3} \textbf{Input:} $S(x,t)$ as well as $s_{i}(x,t)$
	and $e_{i}(x,t)$, for $i=1,2$.
	\begin{enumerate}
		\item Compute the sets $\mathcal{S}_{im}(t)$ of isolated local maxima of
		$\left|s_{i}(\,\cdot\,,t)\right|$ for $i=1,2$. 
		\item Compute attracting OECSs ($e_{2}$--lines) as solutions of the ODE
		\[
		\begin{cases}
		x'(s)=\mathrm{sign}\left\langle e_{2}(x(s)),\tfrac{dx(s-\Delta)}{ds}\right\rangle e_{2}(x(s))\\
		x(0)\in\mathcal{S}_{1m}.
		\end{cases}
		\]
		Stop integration when $\left|s_{1}(x(s))\right|$ ceases to be monotone
		decreasing. 
		\item Compute repelling OECSs ($e_{1}$--lines) as solutions of the ODE
		\[
		\begin{cases}
		x'(s)=\mathrm{sign}\left\langle e_{1}(x(s)),\tfrac{dx(s-\Delta)}{ds}\right\rangle e_{1}(x(s))\\
		x(0)\in\mathcal{S}_{2m}.
		\end{cases}
		\]
		Stop integration when $\left|s_{2}(x(s))\right|$ ceases to be monotone
		decreasing. 
	\end{enumerate}
	\textbf{Output: }Hyperbolic OECSs. 
\end{algorithm}

\subsection{Parabolic OECSs}

\begin{algorithm}[H]
\caption{Compute Parabolic OECSs}
	\label{algorithm4} \textbf{Input:} $S(x,t)$ as well as $s_{i}(x,t)$
	and $e_{i}(x,t)$, for $i=1,2$, and the position and type (wedge
	or trisector) of the rate-of-strain singularities $x^{j}(t)$ satisfying
	$s_{1}(x^{j}(t),t)=s_{2}(x^{j}(t),t),$ $j=1,\ldots,N.$
	\begin{enumerate}
		\item For each trisector-type singularity $x^{j}(t)$, compute the $e_{1}$-lines
		and $e_{2}$-lines connecting the trisector to a wedge. 
		\item Out of these heteroclinic connections, keep only those that are weak
		minimizers of the corresponding neutrality function $\mathcal{N}_{e_{i}}$,
		$i=1,2$. (cf. eqs. (\ref{eq:Neutralitydef}-\ref{eq:WeakminimofNeutrality2})).
		\item With the separatrices obtained in this fashion, build alternating
		chains of heteroclinc $e_{1}$-lines and $e_{2}$-lines
	\end{enumerate}
	\textbf{Output: }Parabolic OECSs. 
\end{algorithm}

\section{Example: Mesoscale OECSs in large-scale ocean data}

\label{NumScheme} We use Algorithms 1-4 from section 6 to locate
OECSs in a two-dimensional ocean surface velocity data set derived
from AVISO satellite altimetry measurements (\url{http://www.aviso.oceanobs.com}).
The domain of interest is the Agulhas leakage in the Southern Ocean
bounded by longitudes $[17^{\circ}W,7^{\circ}E]$, latitudes $[38^{\circ}S,22^{\circ}S]$
and the time slice we selected correspond to $t=24\ \mathrm{November}\ 2006$. 

Under the geostrophic assumptions, the ocean surface height measured
by satellites plays the role of a streamfunction for surface currents.
With $h$ denoting the sea surface height, the velocity field in longitude
latitude coordinates, $[\phi,\ \theta]$, can be expressed as 
\[
\dot{\phi}=-\dfrac{g}{R^{2}f(\theta)\cos\theta}\partial_{\theta}h(\phi,\theta,t),\ \ \ \dot{\theta}=\dfrac{g}{R^{2}f(\theta)\cos\theta}\partial_{\phi}h(\phi,\theta,t),
\]
where $f(\theta):=2\Omega\sin\theta$ denotes the Coriolis parameter,
$g$ the constant of gravity, $R$ the mean radius of the hearth and
$\Omega$ its mean angular velocity. The velocity field is available
at weekly intervals, with a spatial longitude-latitude resolution
of $0.25^{\circ}$. For more detail on the data, see \cite{Beron-VeraDatassh2013}.

\subsection{Elliptic OECSs}

Following Algorithm 1 in section 6.1, we locate singularities of the
rate-of-strain tensor $S(x,t)$, and discard isolated wedge singularities
whose distance to the closest wedge point is larger than the typical
mesoscale distance of $2^{\circ}\approx200km$. The remaining wedge
pairs mark candidate regions for elliptic OECSs.

Along with elliptic OECSs, we will also show the value of the Okubo--Weiss
(OW) parameter
\[
OW(x,t)=s_{2}^{2}(x,t)-\omega^{2}(x,t),
\]
where $\omega(x,t)$ denotes the vorticity. Spatial domains with $OW(x,t)<0$
are frequently used indicator of instantaneous ellipticity in unsteady
fluid flows \cite{Okubo1970,Weiss1991}. While the OW parameter is
not objective (the vorticity term will change under rotations), its
simplicity makes it broadly used in locating coherent vortices. 

\begin{figure}[H]
\subfloat[]{\includegraphics[width=.95\textwidth]{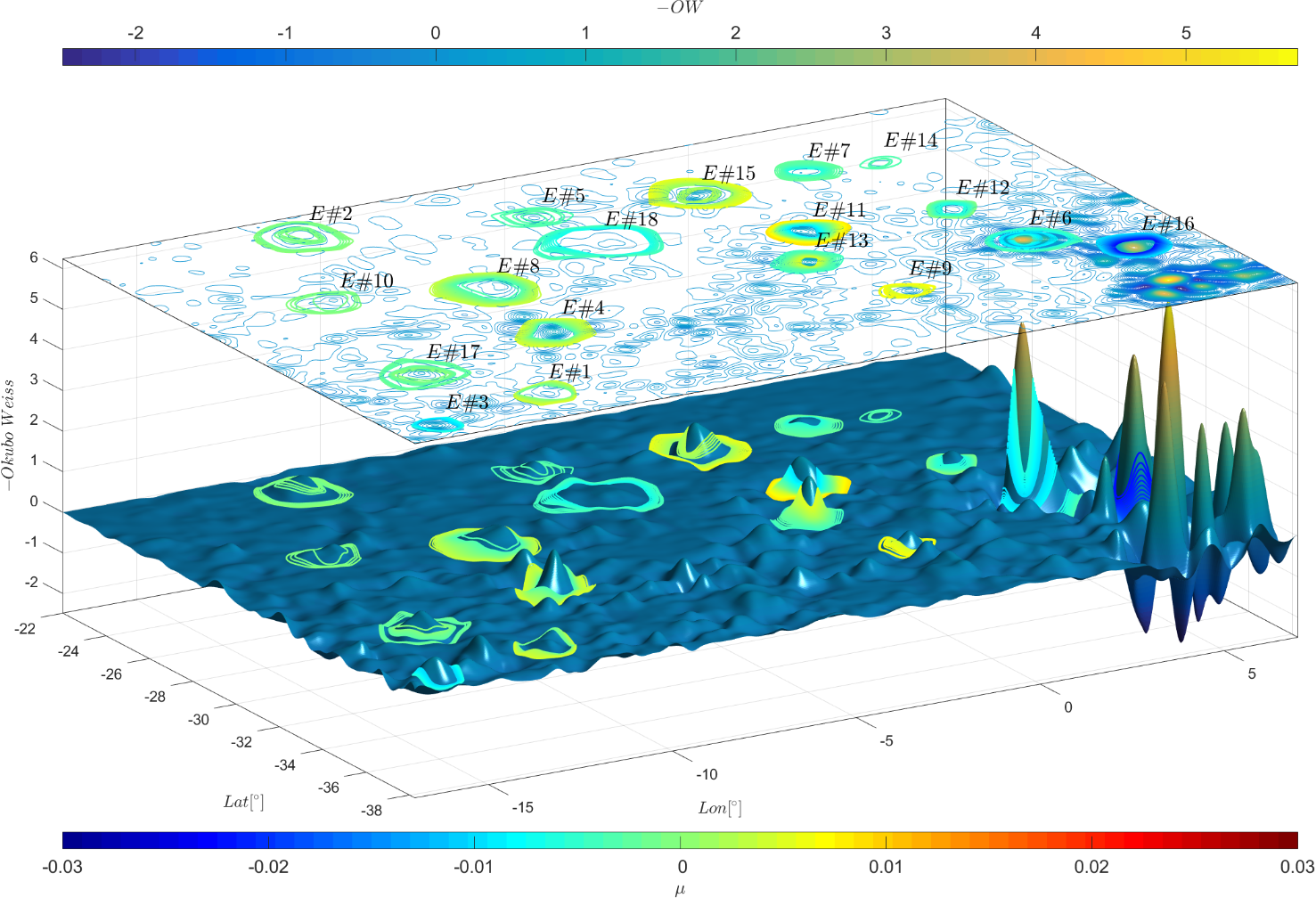}\label{fig:ElliptonOWSurf}}\\
\centering{\subfloat[]{\includegraphics[width=0.35\columnwidth]{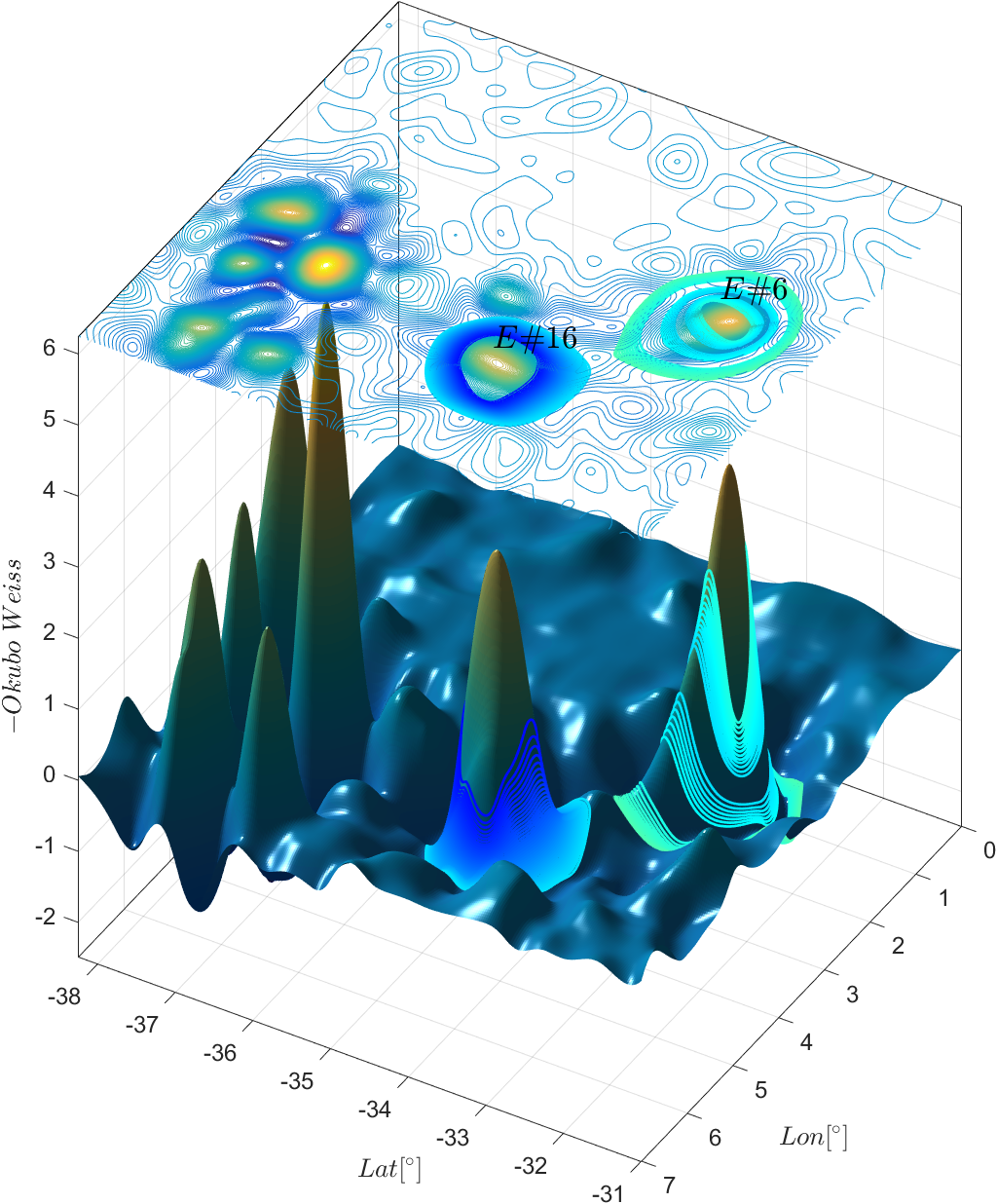}\label{fig:ElliptonOWSurfZoom}}}
	\caption{(a) Elliptic OECSs for different values of stretching rate $\mu$
		(bottom colorbar) on a surface representing the magnitude of $-OW(x,t)$
		(upper colorbar or z-axis). On top, the same elliptic OECSs are shown
		over the level sets of $-OW$. Black numbers label vortical regions
		foliated by families of elliptic OECSs. (b) A different view for the
		bottom-right part of the same domain.}

	\label{fig:ALLISGonOW3dsurf}
\end{figure}

In the domain of study, we obtain a total of eighteen objectively
detected vortical regions, each filled with families of elliptic OECSs.
We plot these families over a two-dimensional graph of $-OW(x,t)$,
and also project them onto the level curves of $-OW(x,t)$ in the
plane (Fig. \ref{fig:ALLISGonOW3dsurf}).

Note that the objective vortical regions $E\#18,E\#8$ and $E\#2$
arise in regions where $-OW(x,t)$ is nearly zero and hence indicates
no significant vortical activity. At the same time, at the bottom
right region of the domain, $OW(x,t)$ attains several strong local
minima, even though there are no elliptic OECSs present (cf. Fig.
\ref{fig:ElliptonOWSurfZoom}). 
\begin{figure}[h]
	\subfloat[]{\includegraphics[height=0.2\paperheight]{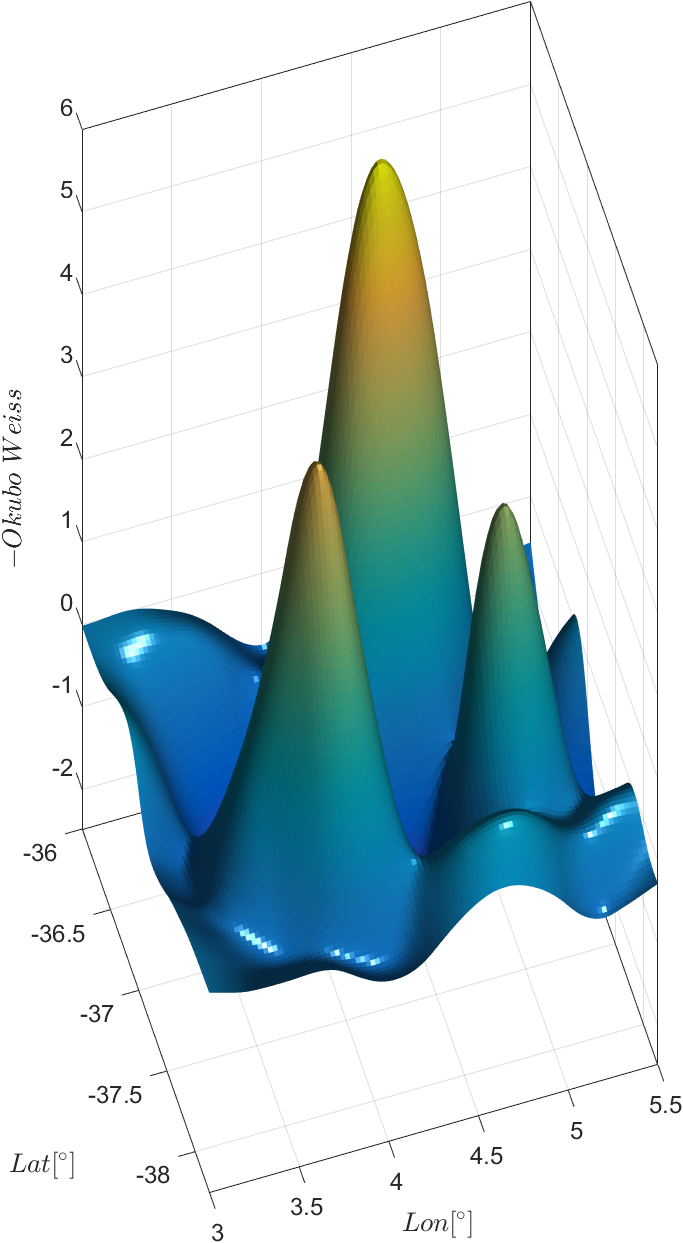}\label{fig:OWBottRight}}\hfill{}
	\subfloat[]{\includegraphics[height=0.18\paperheight]{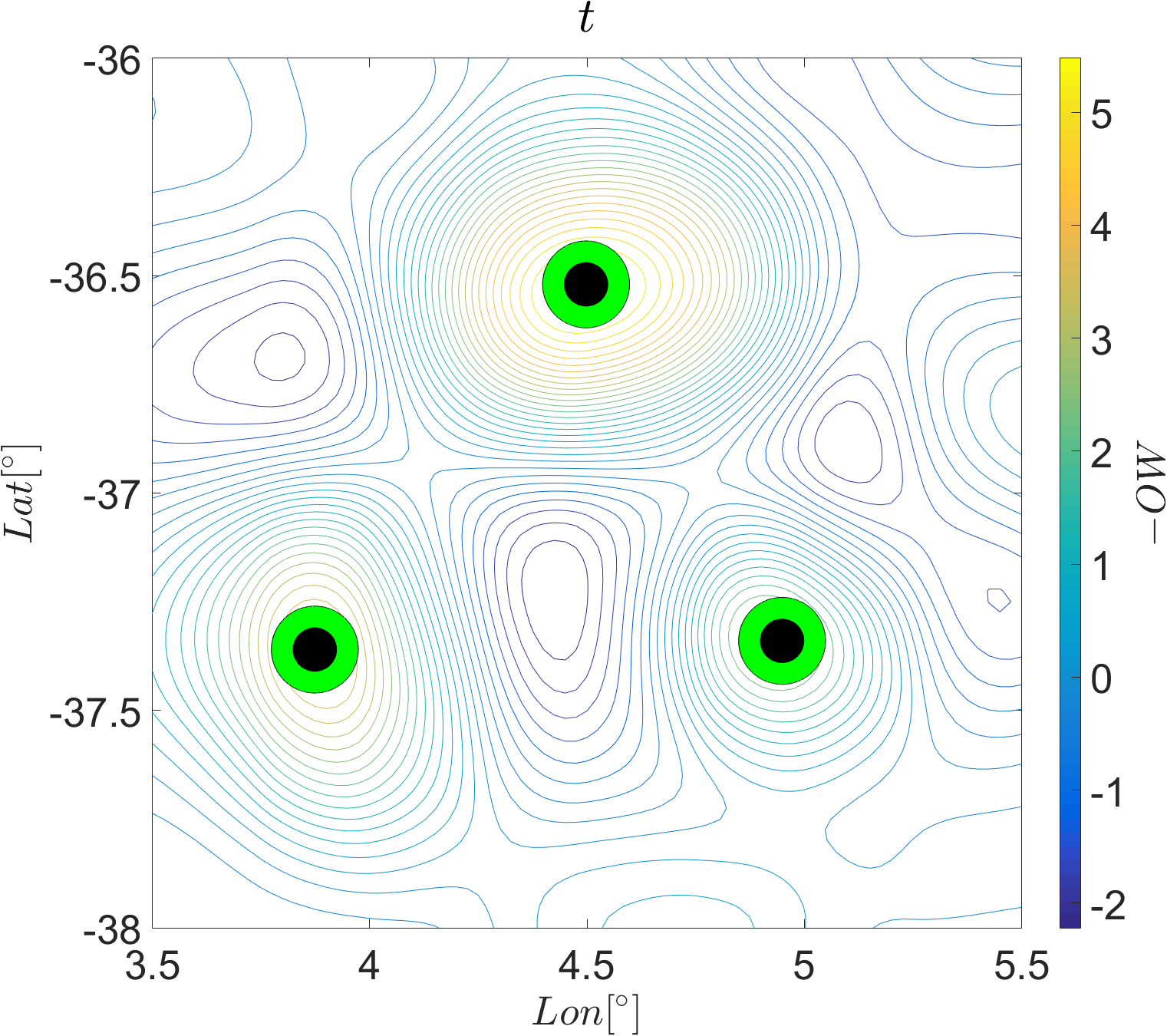}\label{fig:OWBlobst0}}\hfill{}
	\subfloat[]{\includegraphics[height=0.18\paperheight]{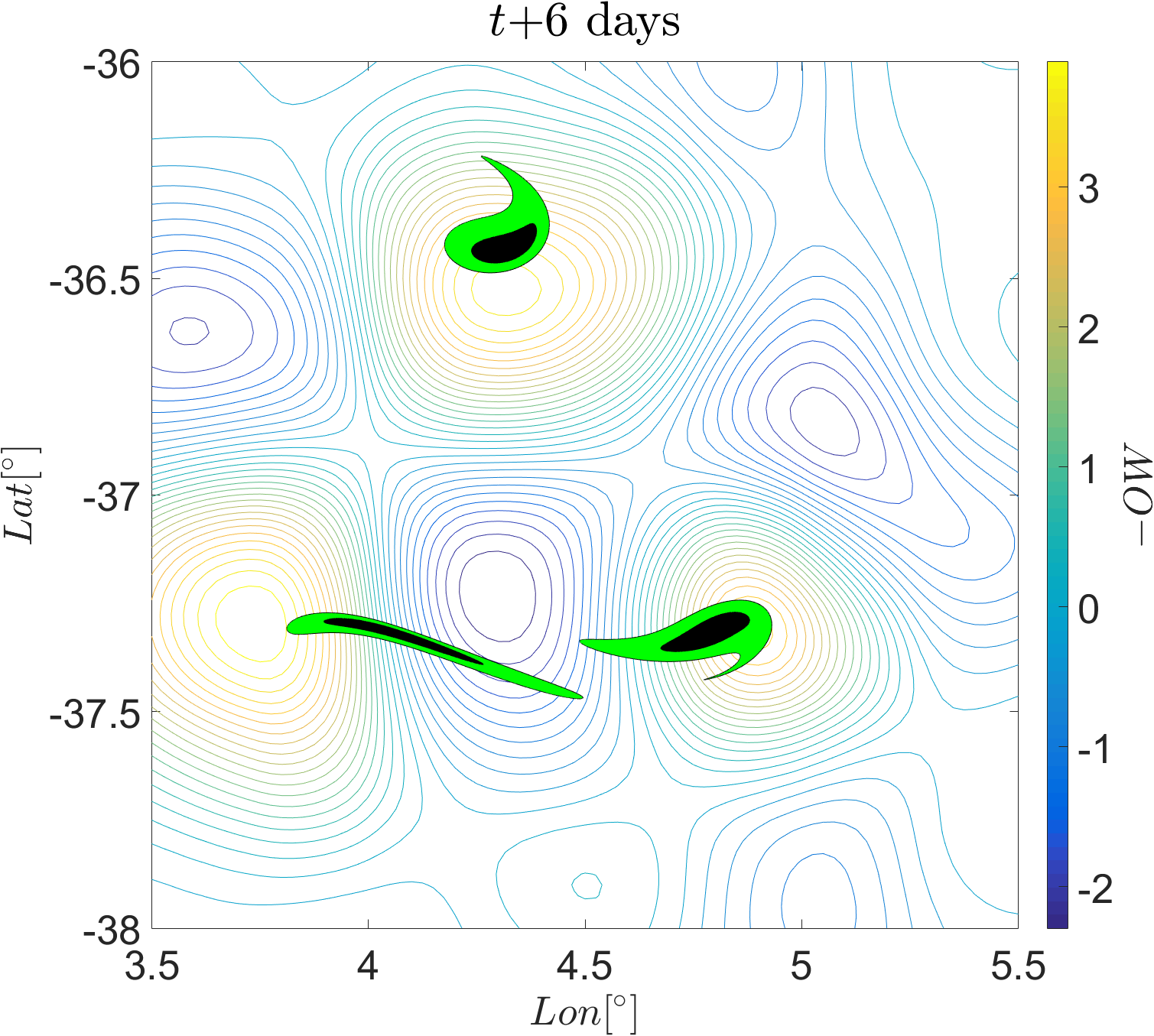}\label{fig:OWBlobstf}}\\	
	\caption{\label{fig:FalsePosOW}(a) Strongest local minima of $OW(x,t)$ in the bottom right region
		of the domain. (b) Material blobs of initial conditions centered on
		local minima of $OW(x,t)$. (c) Deformed material blobs after six
		days.}
\end{figure}
As a representative example, we show
in Fig. \ref{fig:FalsePosOW} the three strongest local minima of
$OW(x,t)$ in this region along with the deformation of material blobs,
initially centered on those minima, for an integration time of six
days. 
For comparison, Fig. \ref{fig:StretchInElliptECSs} shows the deformation
experienced by blobs initialized within two elliptic OECSs $(E\#11,E\#13)$
for different integration times, up to six days.
\begin{figure}[h]
	\subfloat[]{\includegraphics[width=0.20\textwidth]{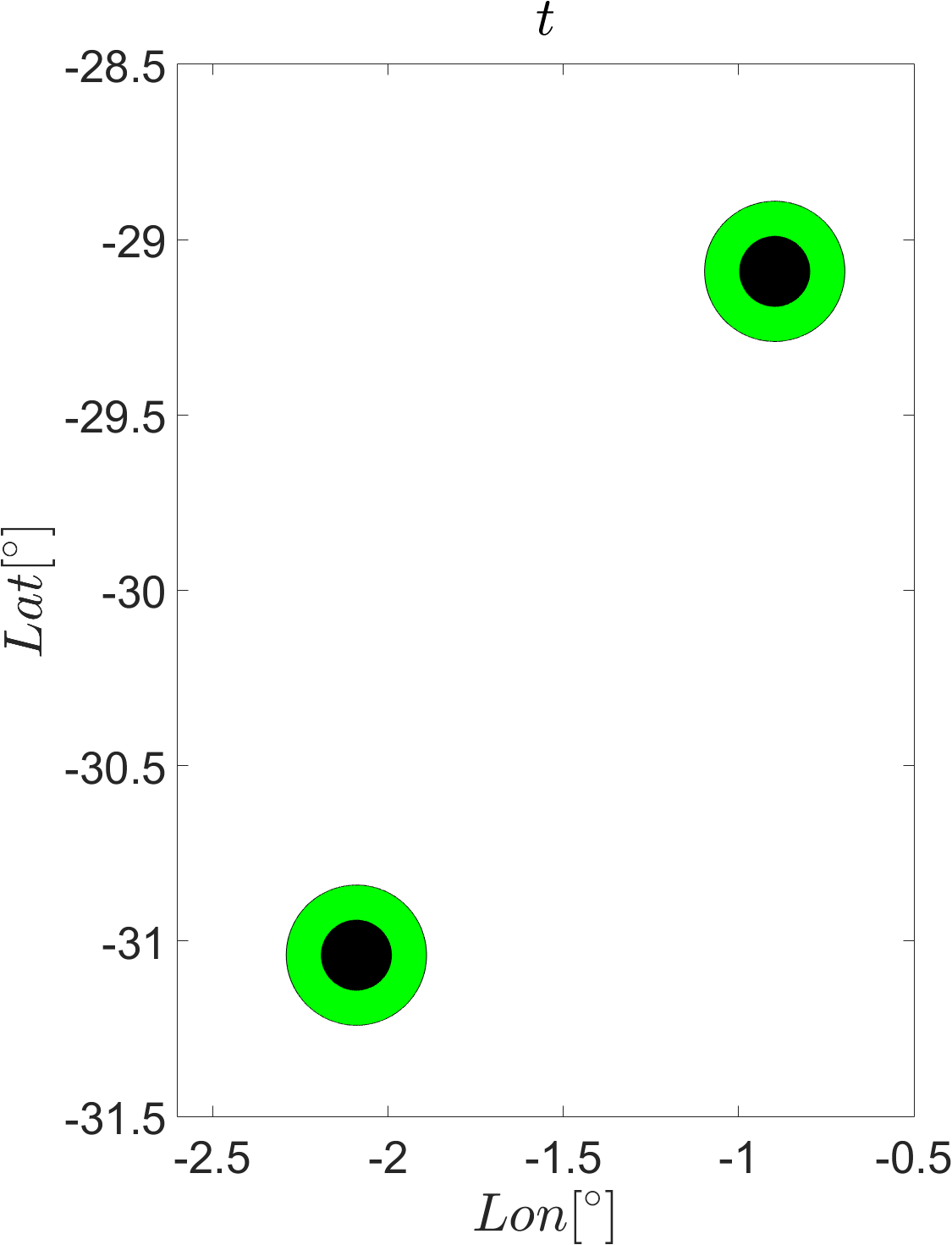}\label{fig:SretchEllip0}}\hfill{}
	\subfloat[]{\includegraphics[width=0.20\textwidth]{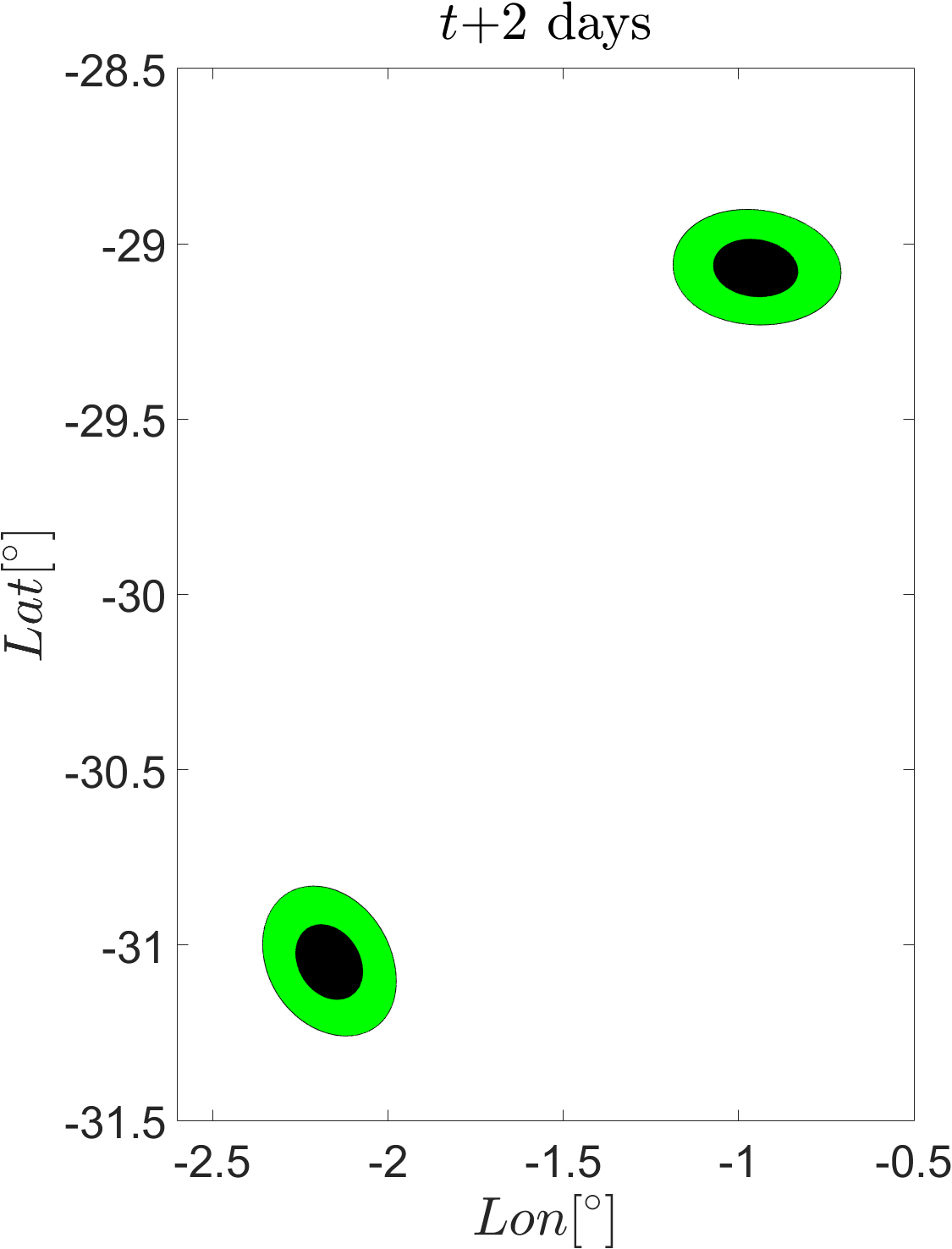}\label{fig:SretchEllip2}}\hfill{}
	\subfloat[]{\includegraphics[width=0.20\textwidth]{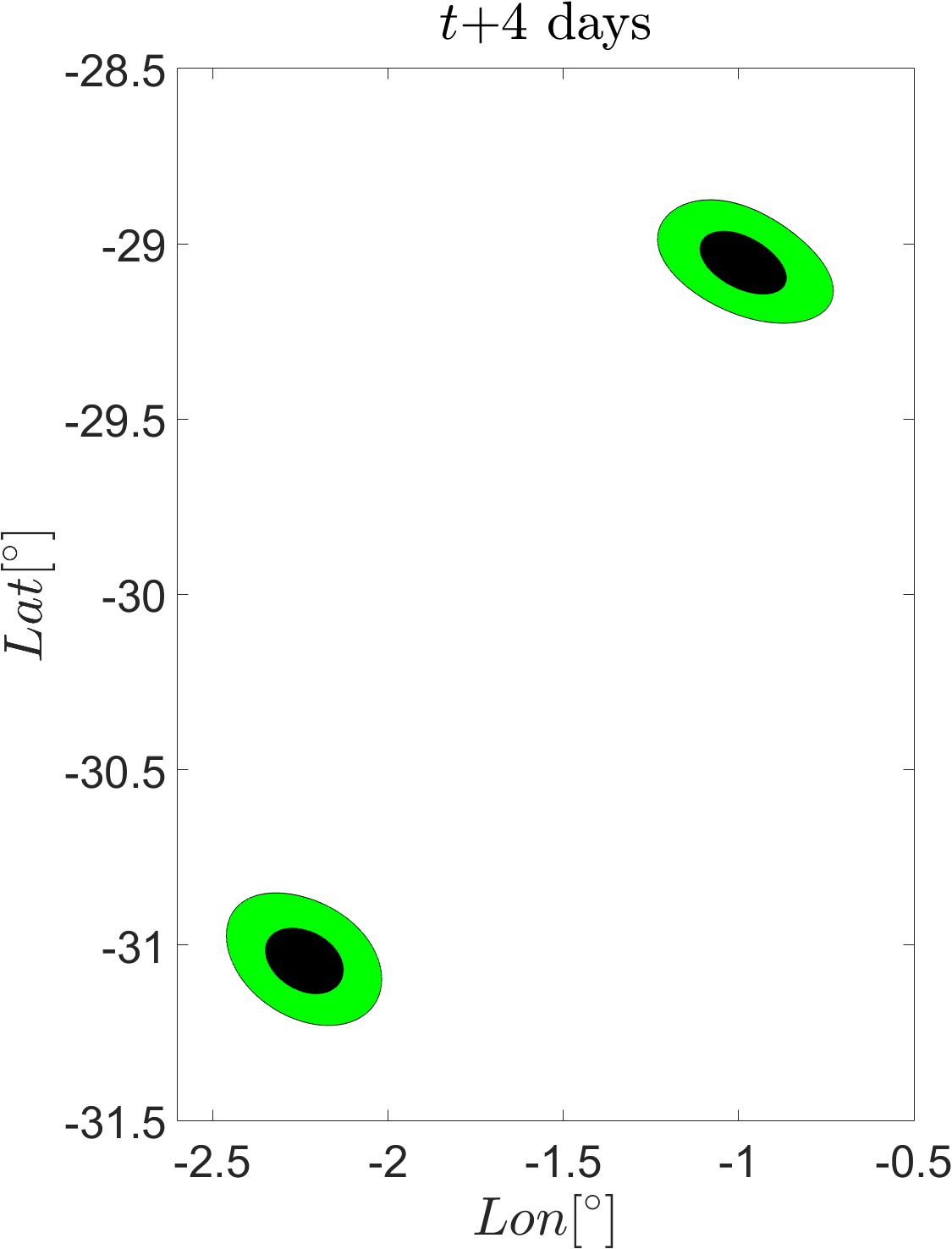}\label{fig:SretchEllip4}}\hfill{}
	\subfloat[]{\includegraphics[width=0.20\textwidth]{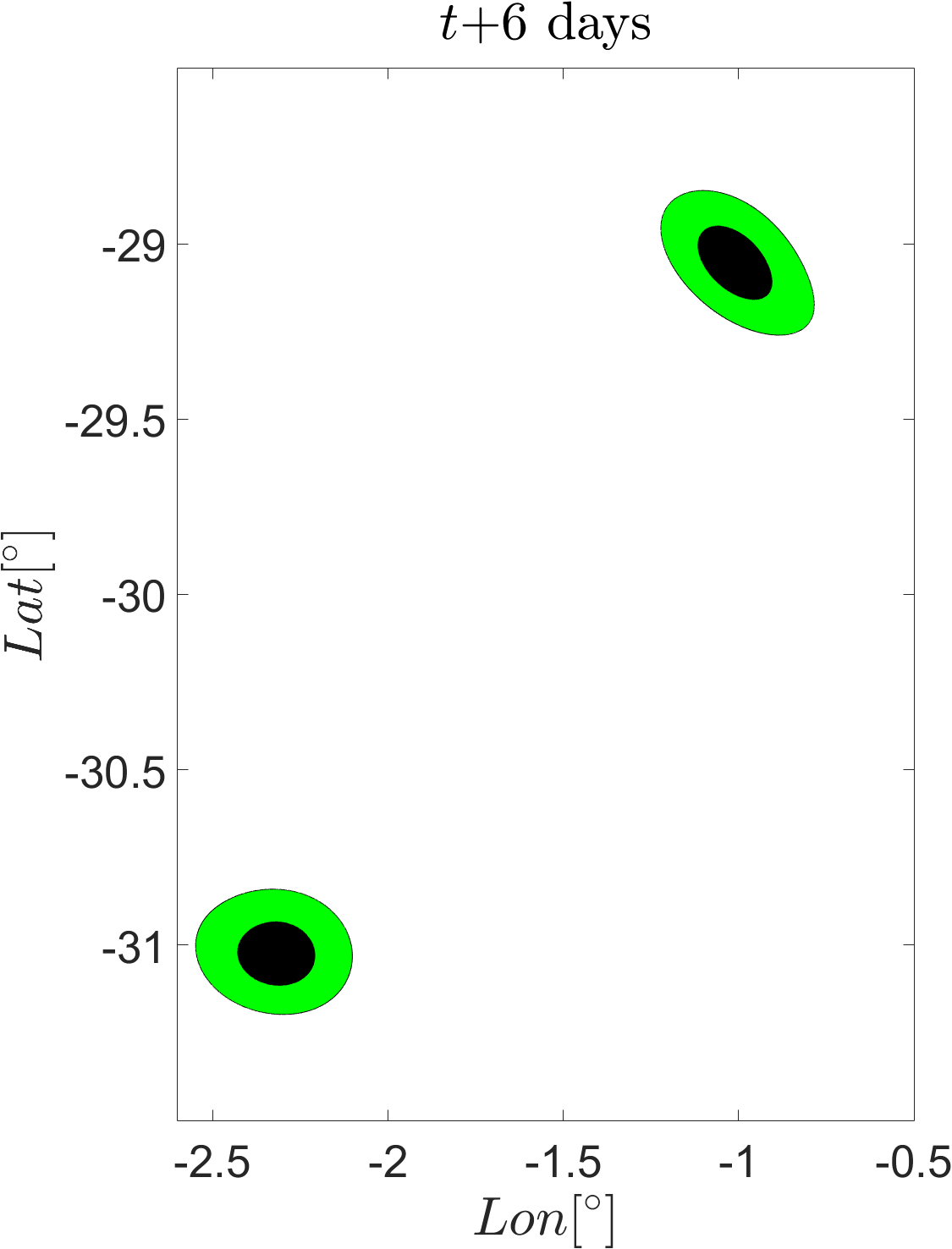}\label{fig:SretchEllip6}}\\
	\caption{Moderate deformation experienced of material blobs of initial conditions
		released within the elliptic OECSs $E\#11$ and $E\#13$.}
	\label{fig:StretchInElliptECSs}
\end{figure}

The blobs in Fig. \ref{fig:StretchInElliptECSs} barely deform, and
hence instantaneously computed elliptic OECSs exhibit short-term Lagrangian
vortex-type behavior, as expected.

Hence, the regions of highest ellipticity for the Okubo-Weiss parameter
show significant stretching, while regions identified by elliptic
OECSs remain highly coherent over the same time interval. We have,
therefore, clear examples of both false positives and false negatives
for Okubo--Weiss-based vortex detection.

The vortical regions $\#2,\#6,\#8,\#11,\#13,\#15,\#18$ approximate
locations where exceptionally coherent Lagrangian coherent eddies
have been found in other studies (see \cite{AutoMdetectFlorian},
\cite{BlackHoleHaller2013}). These Lagrangian studies cover a time
interval of three months, with their initial time coinciding with
the time of the present Eulerian analysis. Remarkably, about one third
of the elliptic OECSs we find are signatures of elliptic LCSs with
long-term coherence. 

\begin{figure}[h]
\includegraphics[width=1\textwidth]{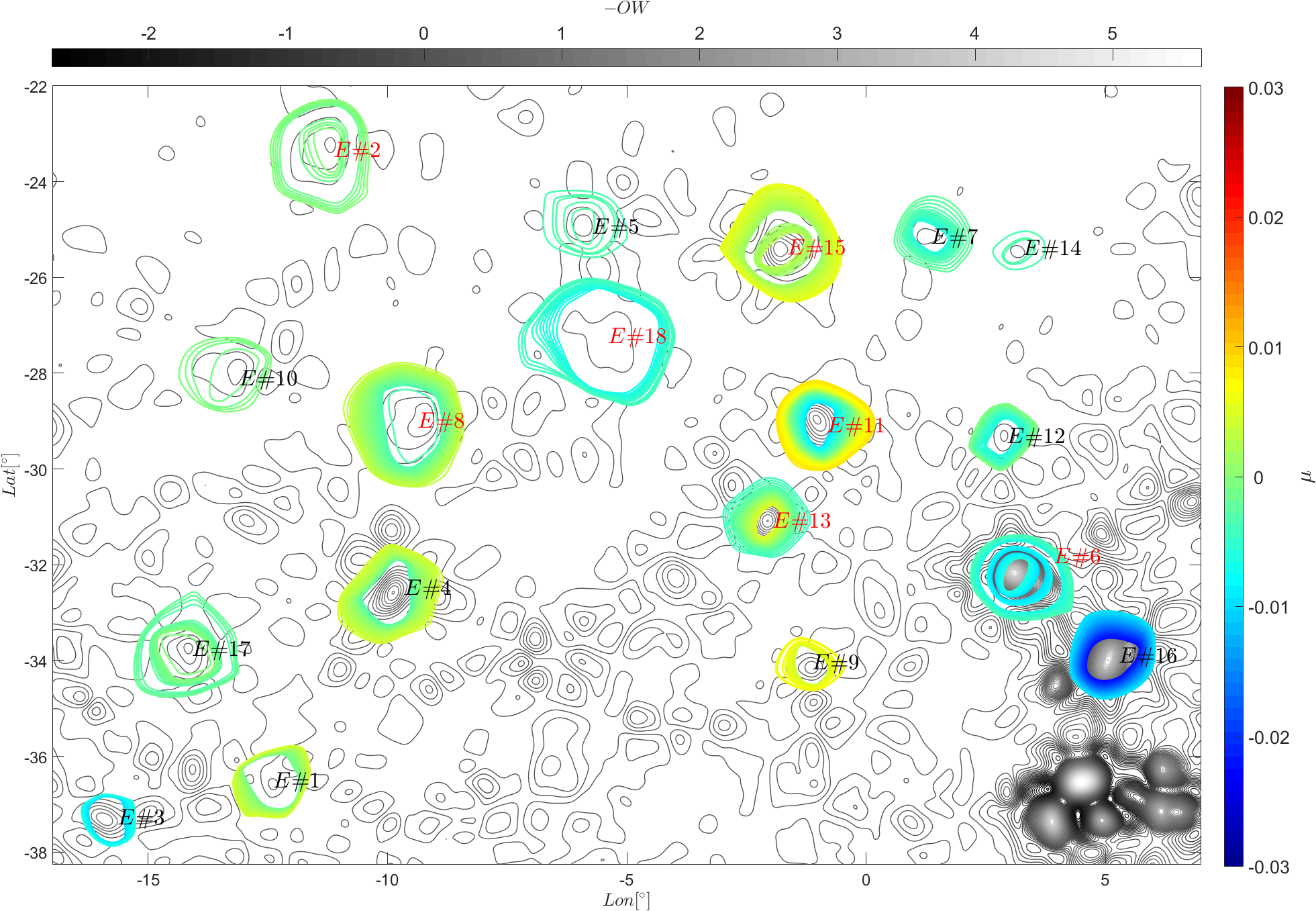}
	\caption{Top view of Fig. \ref{fig:ElliptonOWSurf} with a different colormap
		for the\textbf{ }negative $OW(x,t)$ parameter. Red numbers identify
		elliptic OECSs located in regions where exceptionally coherent Lagrangian
		eddies have been found. }

	\label{fig:ALLISGonOW3dsurfTopview}
\end{figure}

Figure \ref{fig:ALLISGonOW3dsurfTopview} underlines this observation
by showing the top view of Fig. \ref{fig:ElliptonOWSurf} separately,
using a different colormap for the $-OW(x,t)$ scalar field, and tagging
with red numbers elliptic OECSs in the coherent Lagrangian eddy domains
identified by \cite{BlackHoleHaller2013}. Without exception, these
regions are marked by near-zero values of the stretching rate, indicating
a high degree of Eulerian coherence for the elliptic OECSs.

In contrast, most of these highly coherent Lagrangian eddies have
very moderate signatures in the contour plot of the Okubo--Weiss parameter.
Taking all regions of closed $OW(x,t)$ level sets with comparable
values as predictions for coherent eddies would results in an order-of-magnitude
over-prediction for vortical regions. This is consistent with the
findings of \cite{Beron-VeraDatassh2013}, which reports a roughly
tenfold over-prediction of the actual number of coherent Agulhas eddies
by the Okubo--Weiss criterion.

\subsection{Hyperbolic OECSs}

\begin{figure}[h]
	\subfloat[]{\includegraphics[width=0.47\textwidth]{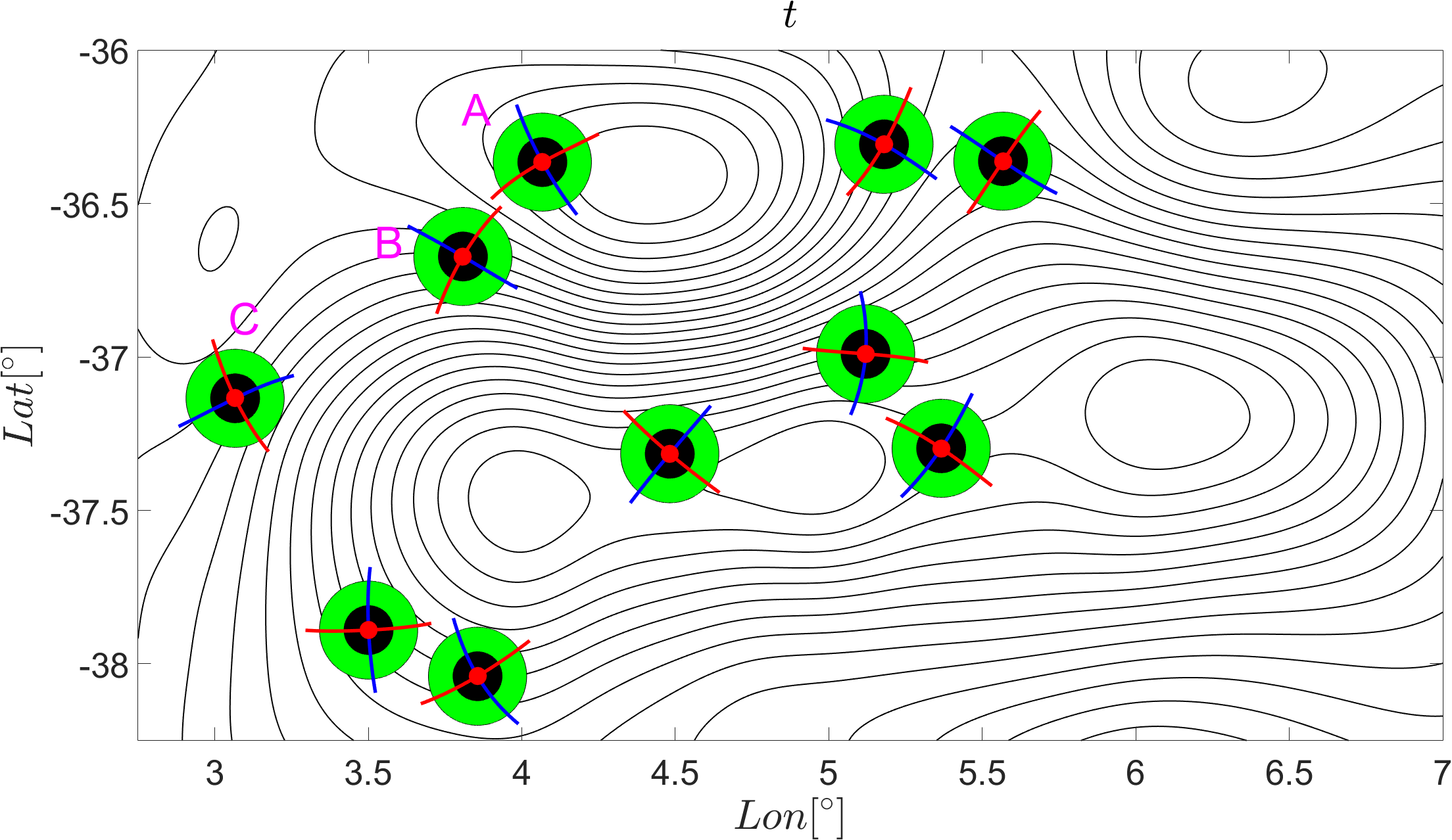}\label{fig:HyplocStrml0}}\hfill{}
	\subfloat[]{\includegraphics[width=0.47\textwidth]{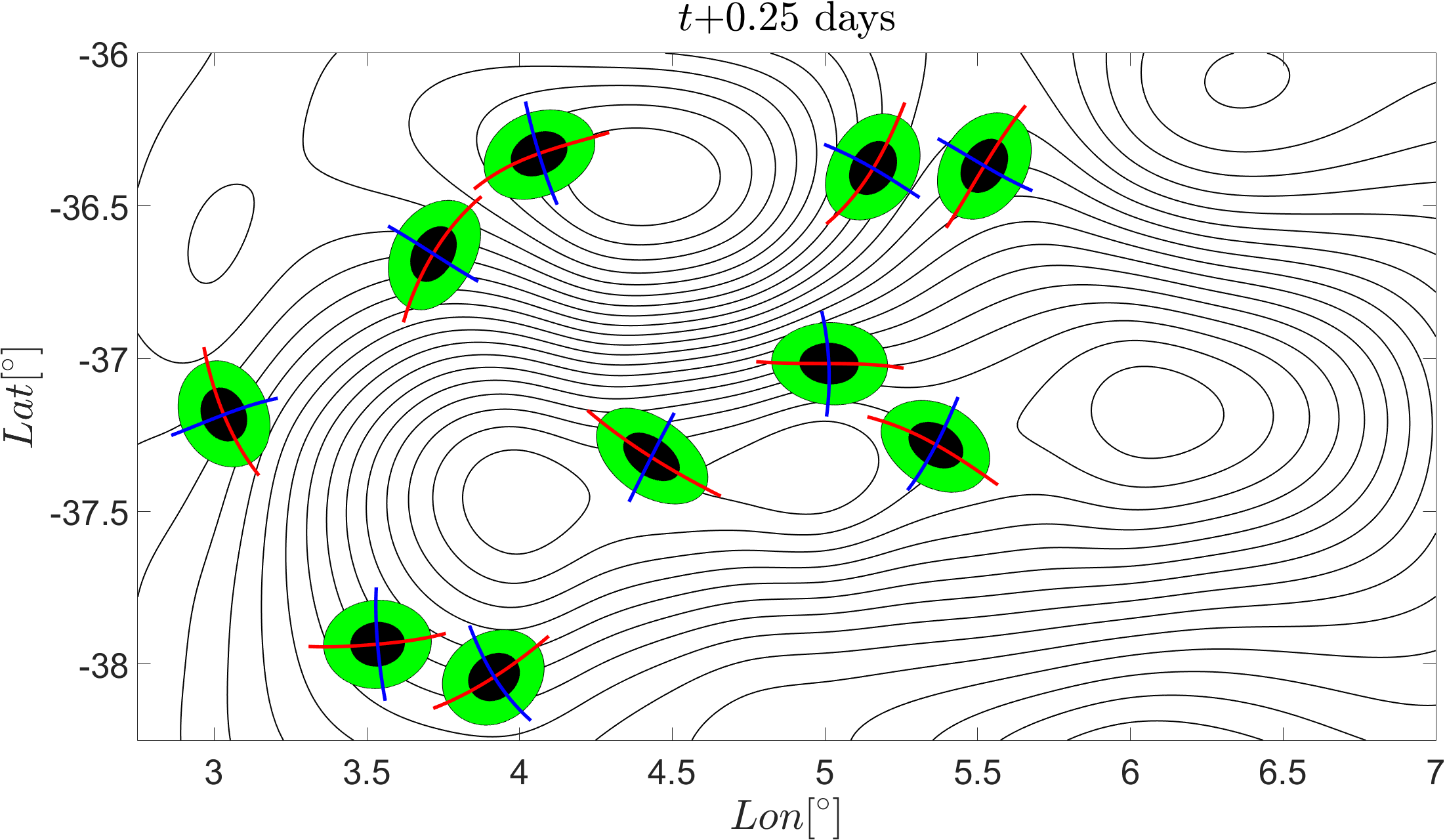}\label{fig:HyplocStrml025}}\\
	\subfloat[]{\includegraphics[width=0.47\textwidth]{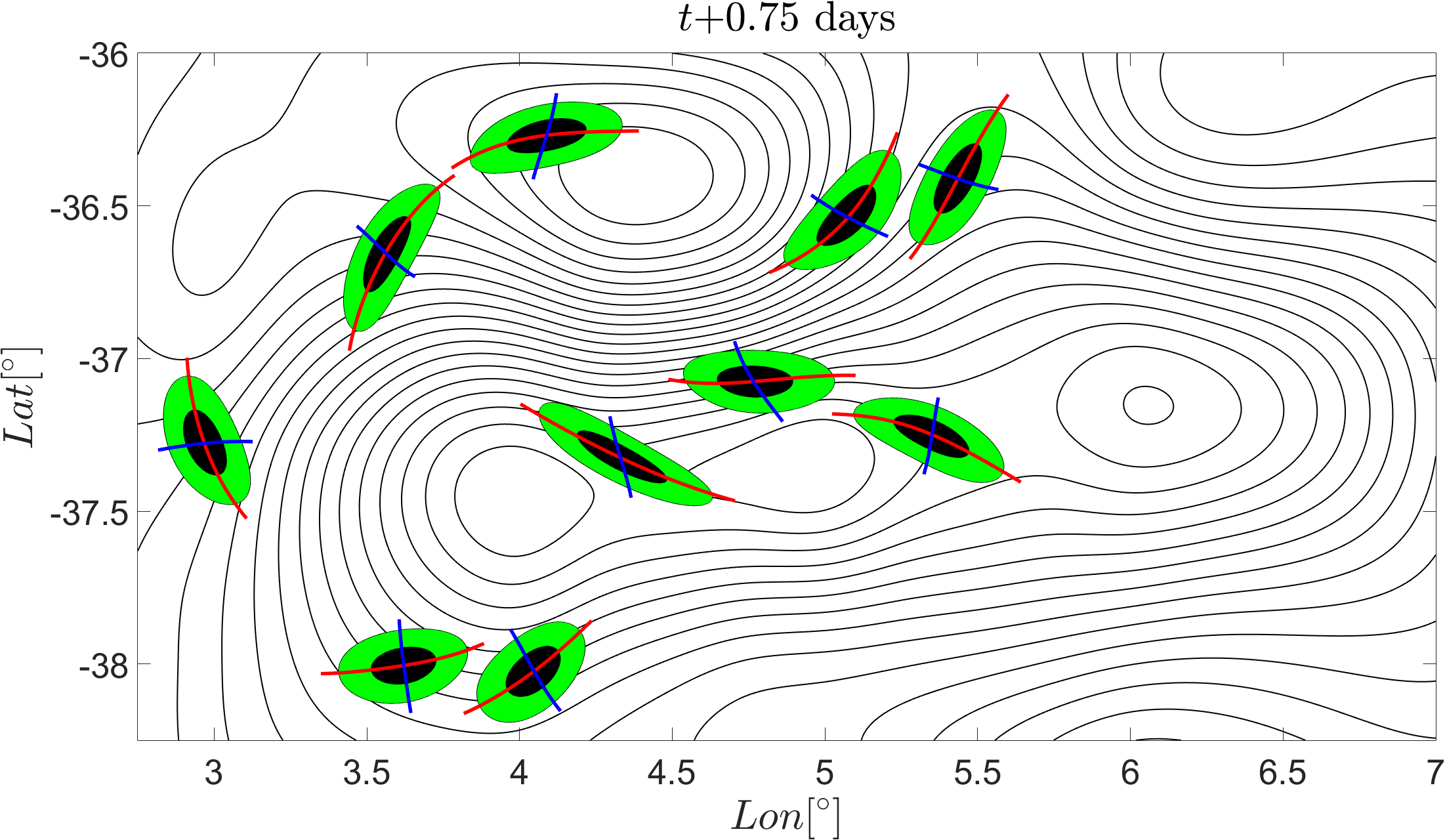}\label{fig:HyplocStrml075}}\hfill{}
	\subfloat[]{\includegraphics[width=0.47\textwidth]{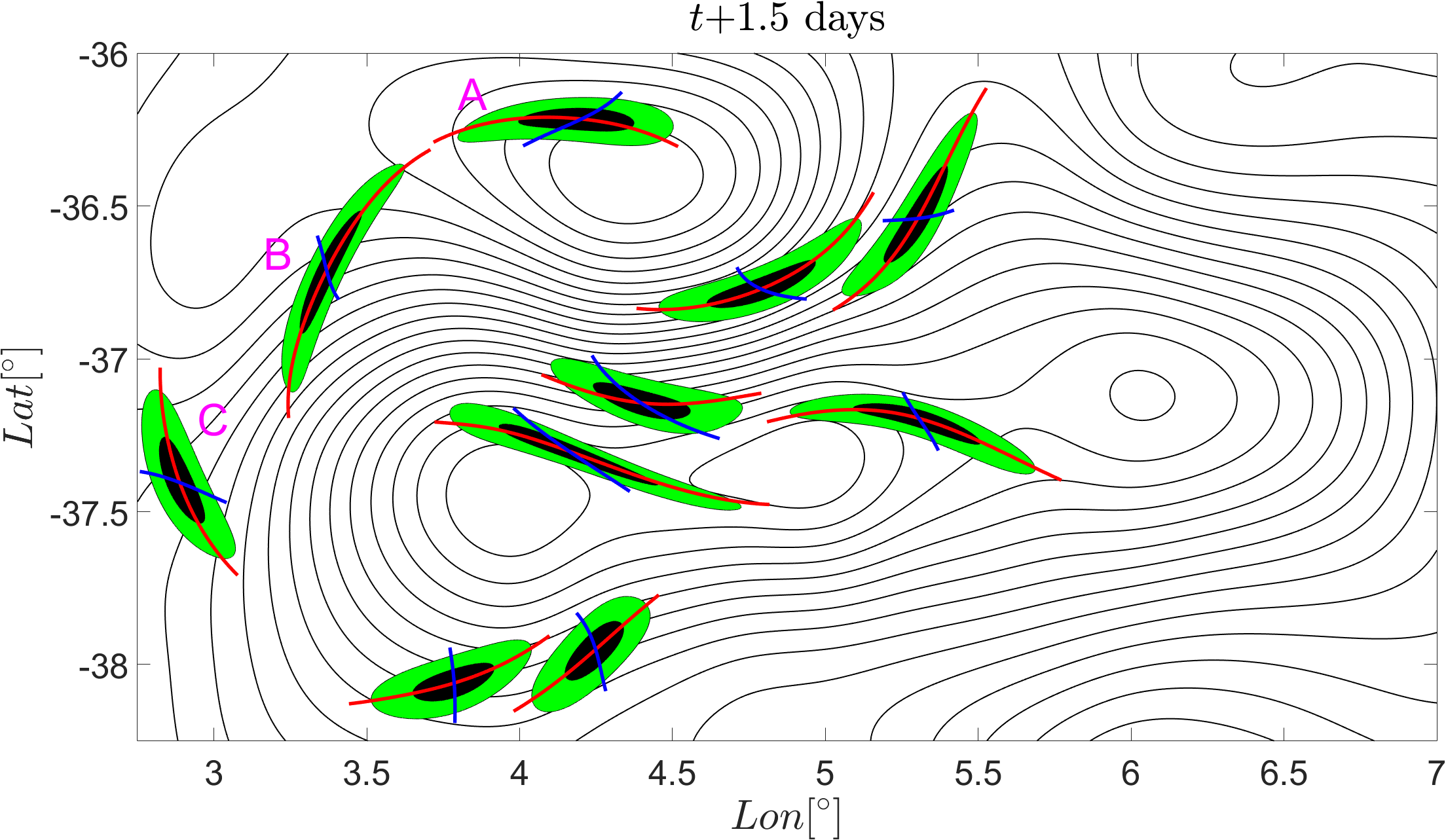}\label{fig:HyplocStrml125}}\\
		
	\caption{(a) Hyperbolic OECSs: attracting (red) and repelling (blue) overlaid
		on streamlines. The red dots denote the\textit{ objective saddle points}
		(cores of the OECSs), i.e., maxima of $s_{2}$, that coincide with
		the minima of $s_{1}$ in the present incompressible flow. (b-d) Advected
		images of the hyperbolic OECSs up to 1.5 days and their effect on
		nearby particles. Magenta letters identify hyperbolic OECSs inducing
		short-term Lagrangian stretching completely hidden in the streamline
		geometry. A: observed hyperbolic behavior within closed instantaneous
		streamlines. B,C: hyperbolic OECSs showing significant stretching,
		with attracting OECSs perpendicular to the streamlines.}
	\label{fig:Hypei12ECSs}
\end{figure}

We compute hyperbolic OECSs at the same time used in our elliptic
OECSs calculations, on a subdomain bounded by longitudes $[2.5^{\circ}E,7^{\circ}E]$
and latitudes $[38.25^{\circ}S,36^{\circ}S]$. Since the velocity
field is incompressible, the maxima of $s_{2}$ coincide with the
minima of $s_{1}$. Consequently, the same OECS cores arise as starting
points in the computation of both attracting and repelling hyperbolic
OECSs. We show all hyperbolic OECSs so obtained in Fig. \ref{fig:Hypei12ECSs}.
As noted earlier, the cores of hyperbolic OECSs represent the \textit{objective
	saddle points }in the unsteady velocity field.\\

Figure \ref{fig:Hypei12ECSs} shows how hyperbolic OECSs act as instantaneous
stable and unstable manifolds for short-term particle motion. Remarkably,
several hyperbolic OECSs cross the local streamlines at large angles
at the initial time (Fig. \ref{fig:HyplocStrml0}), as well as at
later times (Figs. \ref{fig:HyplocStrml025}-\ref{fig:HyplocStrml125}),
explaining the deformation of nearby blobs of fluid. In particular,
the \textit{\emph{objective saddle point}}\emph{ (cf. }\textsl{\emph{Remark
		1}}\emph{)}\textbf{ }captured by a hyperbolic OECS (point A in Figs.
\ref{fig:HyplocStrml0}, \ref{fig:HyplocStrml125}) induces markedly
hyperbolic short-term Lagrangian behavior even though it is located
within an area of closed streamlines. In a similarly surprising fashion,
the \textit{\emph{objective saddle points}} B and C create significant
short-term material stretching in a direction perpendicular to the
local streamlines. 

\begin{figure}[h]
	\centering
\includegraphics[width=0.5\textwidth]{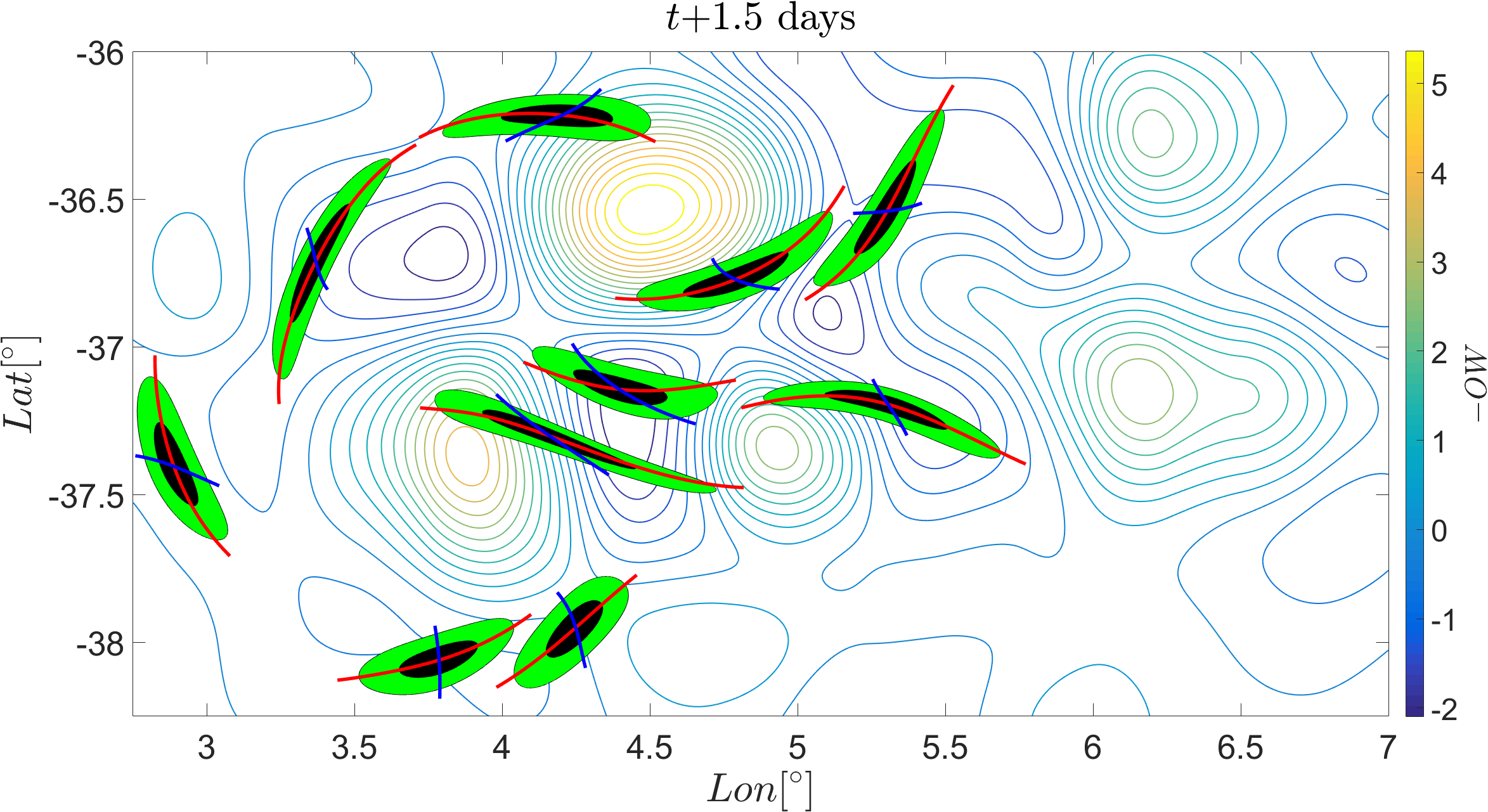}
	\caption{(a) Same figure as Fig. \ref{fig:HyplocStrml125} with level sets
		of $-OW(x,t)$ in the background.}
	\label{fig:Hypei12ECSsOWStrml}
\end{figure}

Figure \ref{fig:Hypei12ECSsOWStrml} presents the same results as
Fig. \ref{fig:HyplocStrml125} but over level sets of the the negative
Okubo--Weiss parameter $-OW(x,t)$. Note how some hyperbolic OECSs
cross closed contours around local minima of $OW(x,t)$ that are generally
believed to signal elliptic regions.

\begin{figure}[h]
	\subfloat[]{\includegraphics[width=0.47\columnwidth]{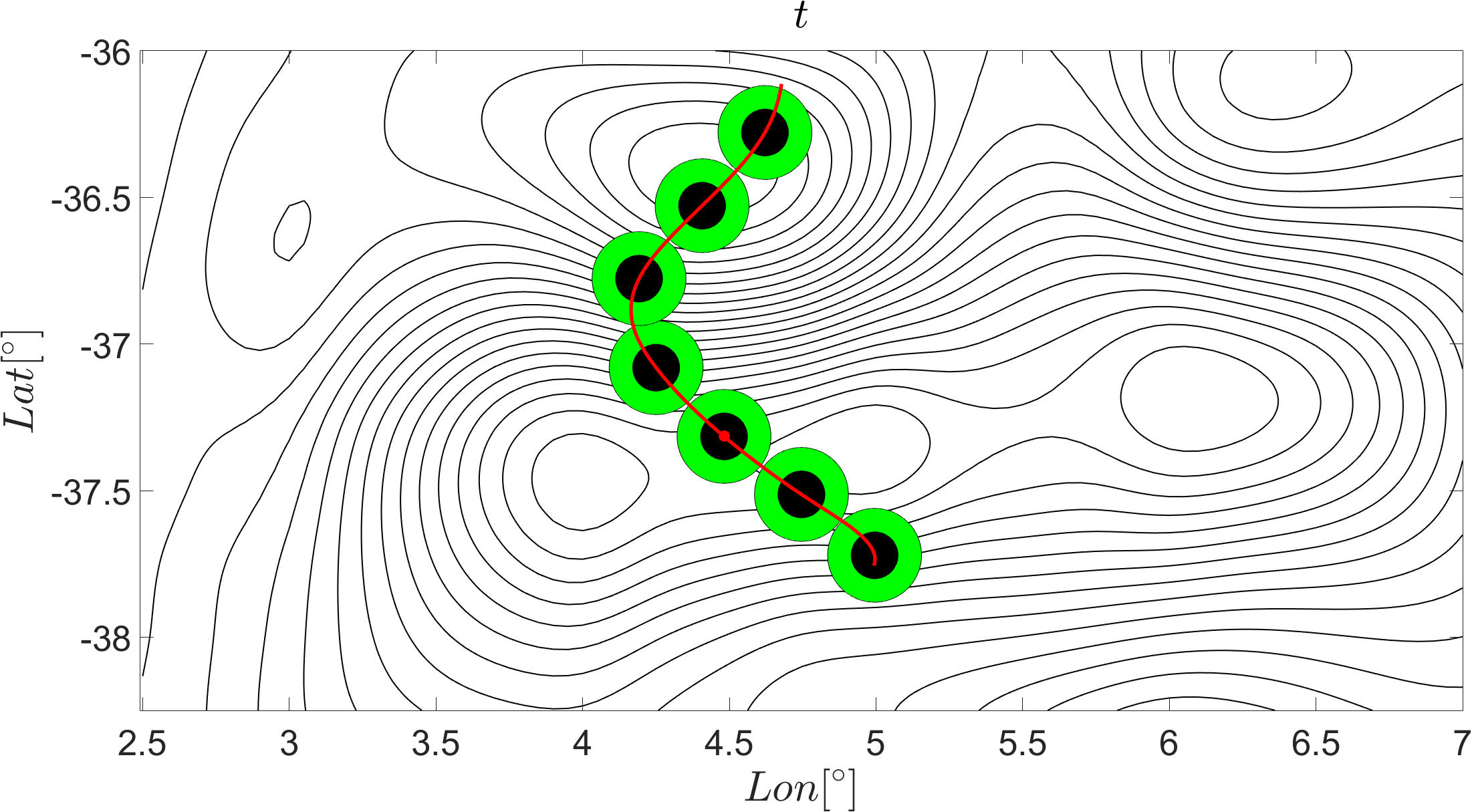}\label{fig:Globxi2t0}}\hfill{}
	\subfloat[]{\includegraphics[width=0.47\columnwidth]{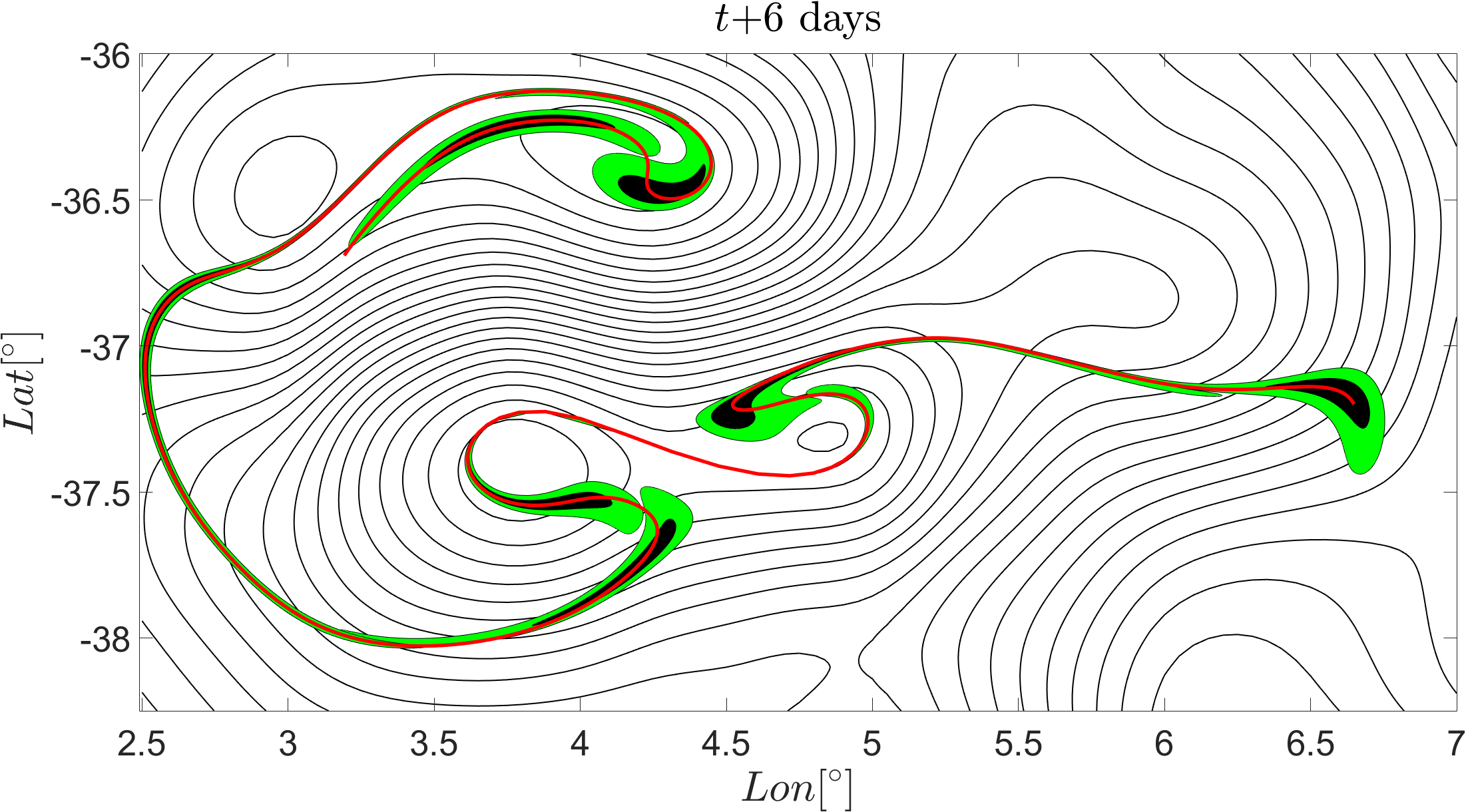}\label{fig:Globxi2t6}}\\	
	\caption{(a) More extended view of an attracting OECS, with blobs of initial
		conditions placed along the OECS, overlaid on streamlines. (b) Advected
		image of the OECS and of the marked initial conditions for 6 days.}
	\label{fig:HypAttractingECSsGlobal}
\end{figure}

Figure \ref{fig:HypAttractingECSsGlobal} illustrates that attracting
OECSs continue to shape short-term tracer deformation patterns in
larger distances from their cores, for times up to six days. Over
this time interval, material blobs align with materially advected
$e_{2}$-lines, underlining the role of attracting OECSs as short-term
unstable manifolds. 

\begin{figure}[h]
	\subfloat[]{\includegraphics[width=0.45\textwidth]{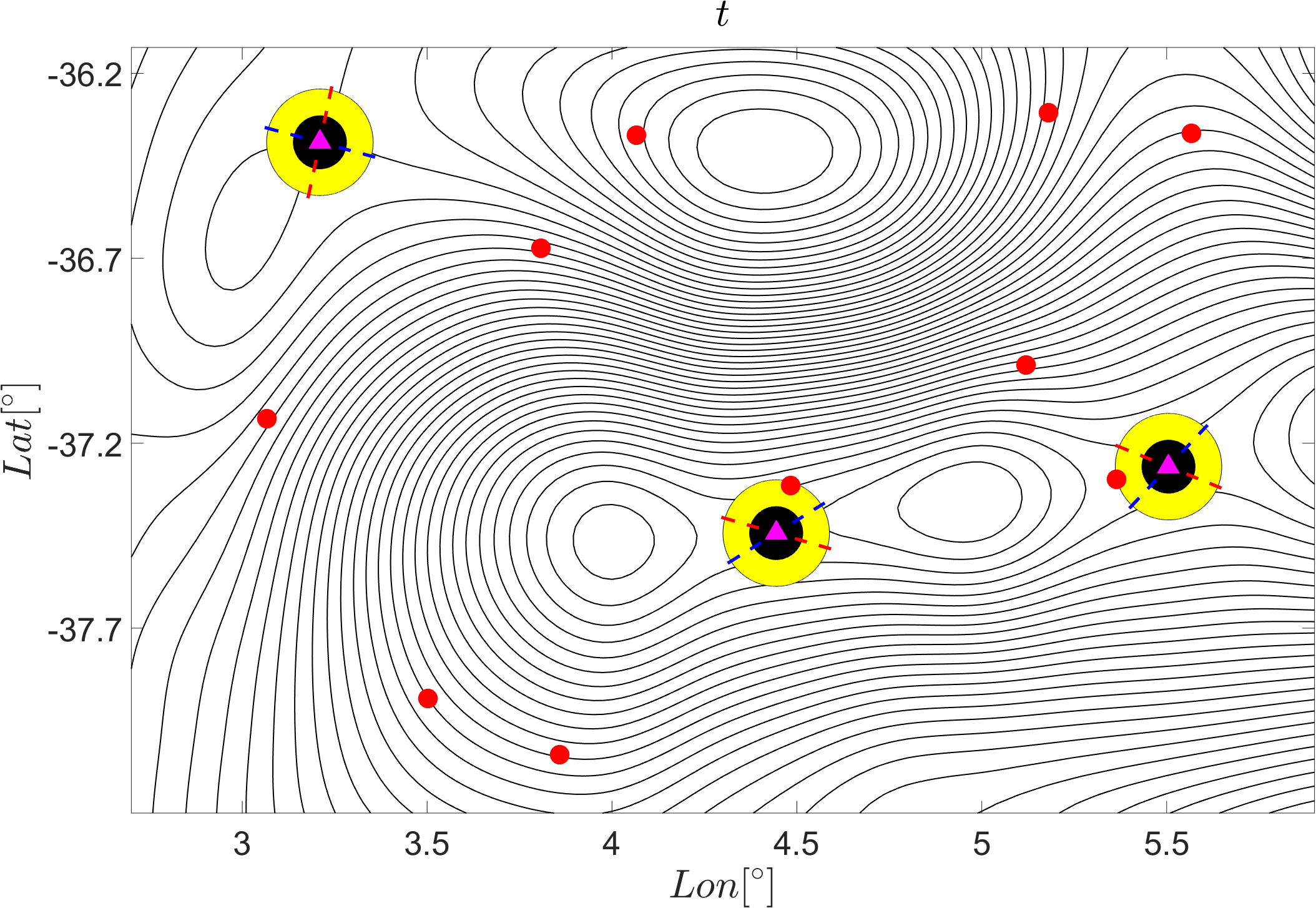}\label{fig:SaddHypStrml0}}\hfill
	\subfloat[]{\includegraphics[width=0.45\textwidth]{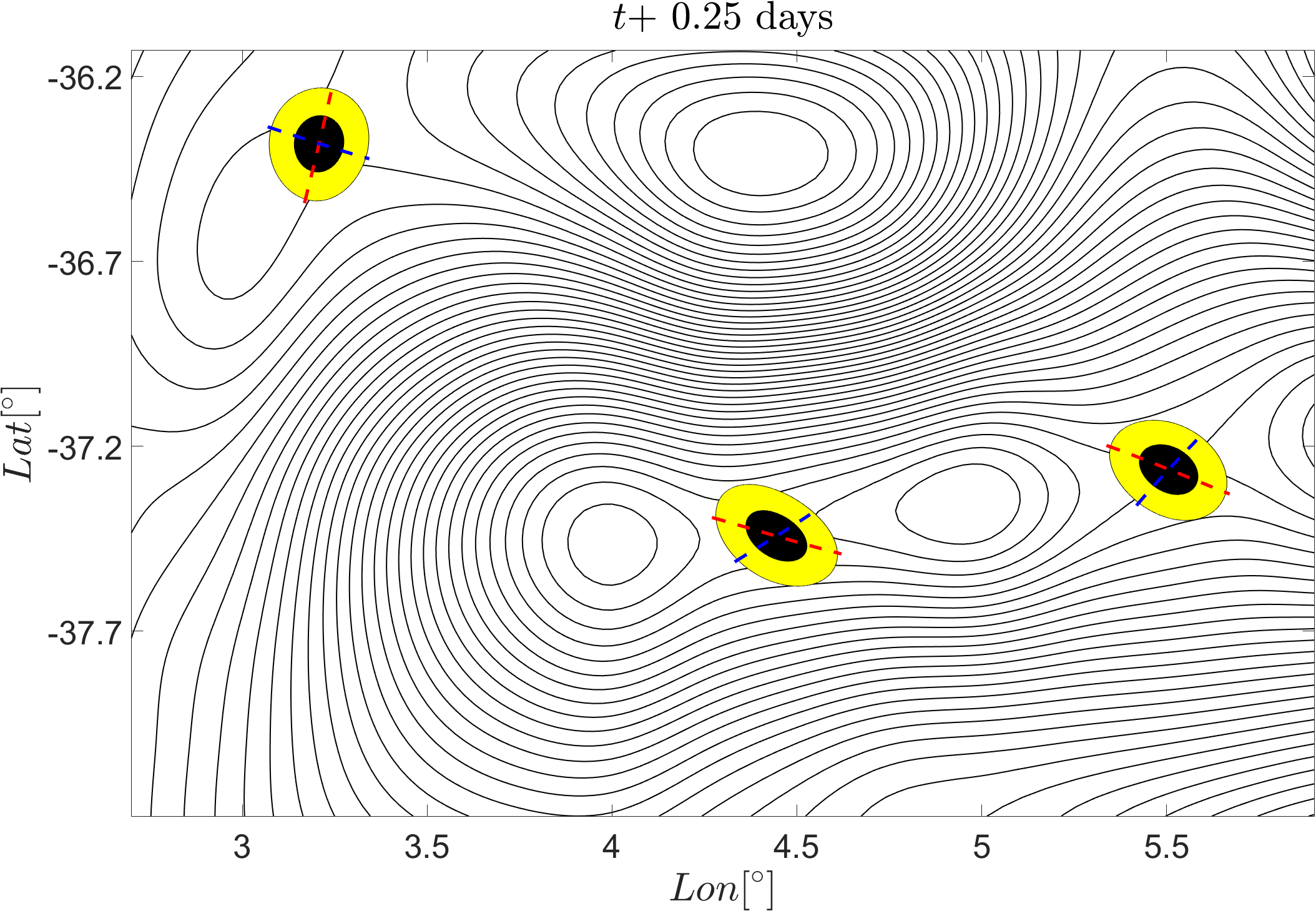}\label{fig:SaddHypStrml025}}\\	
	\subfloat[]{\includegraphics[width=0.45\textwidth]{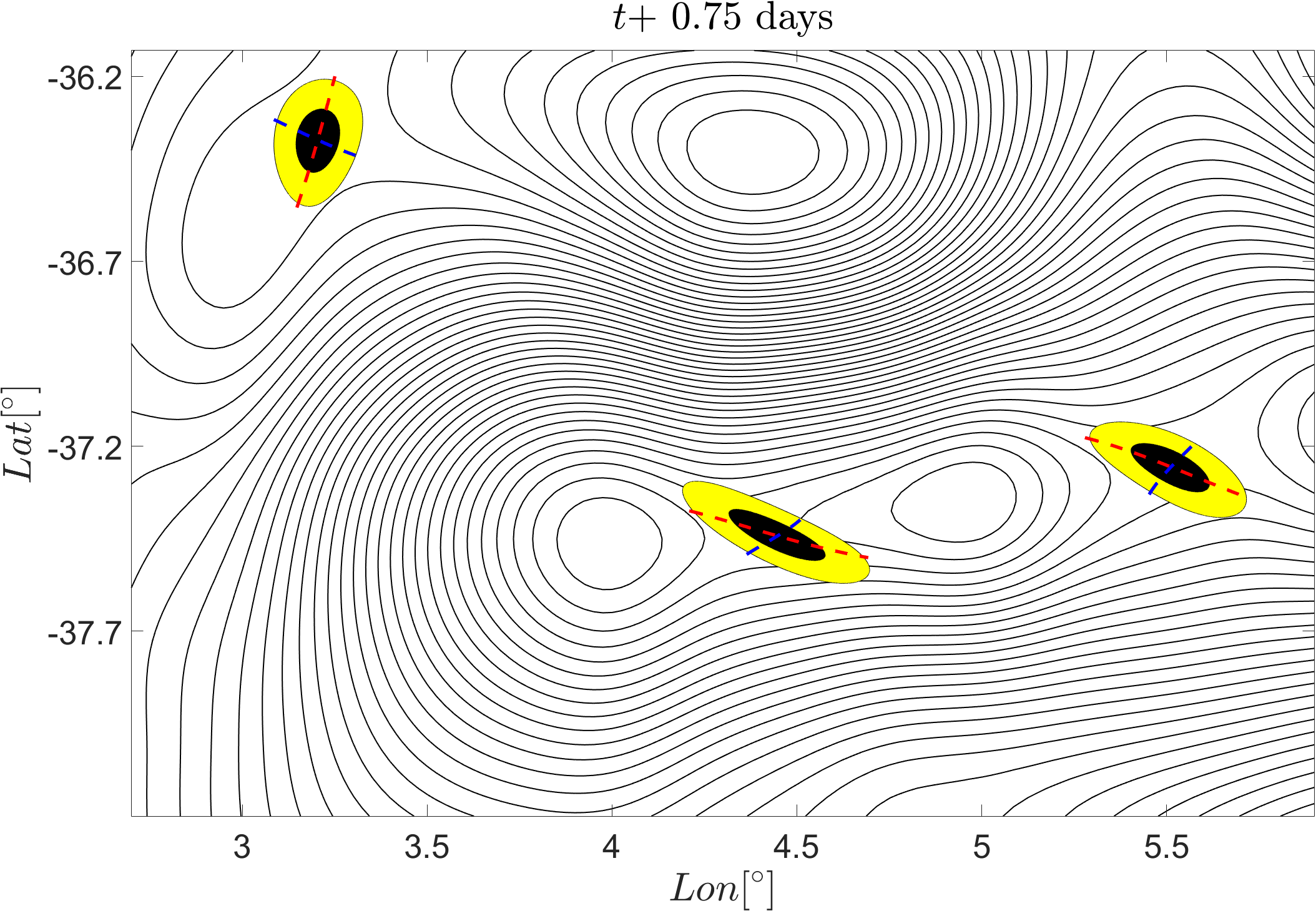}\label{fig:SaddHypStrml075}}\hfill
	\subfloat[]{\includegraphics[width=0.45\textwidth]{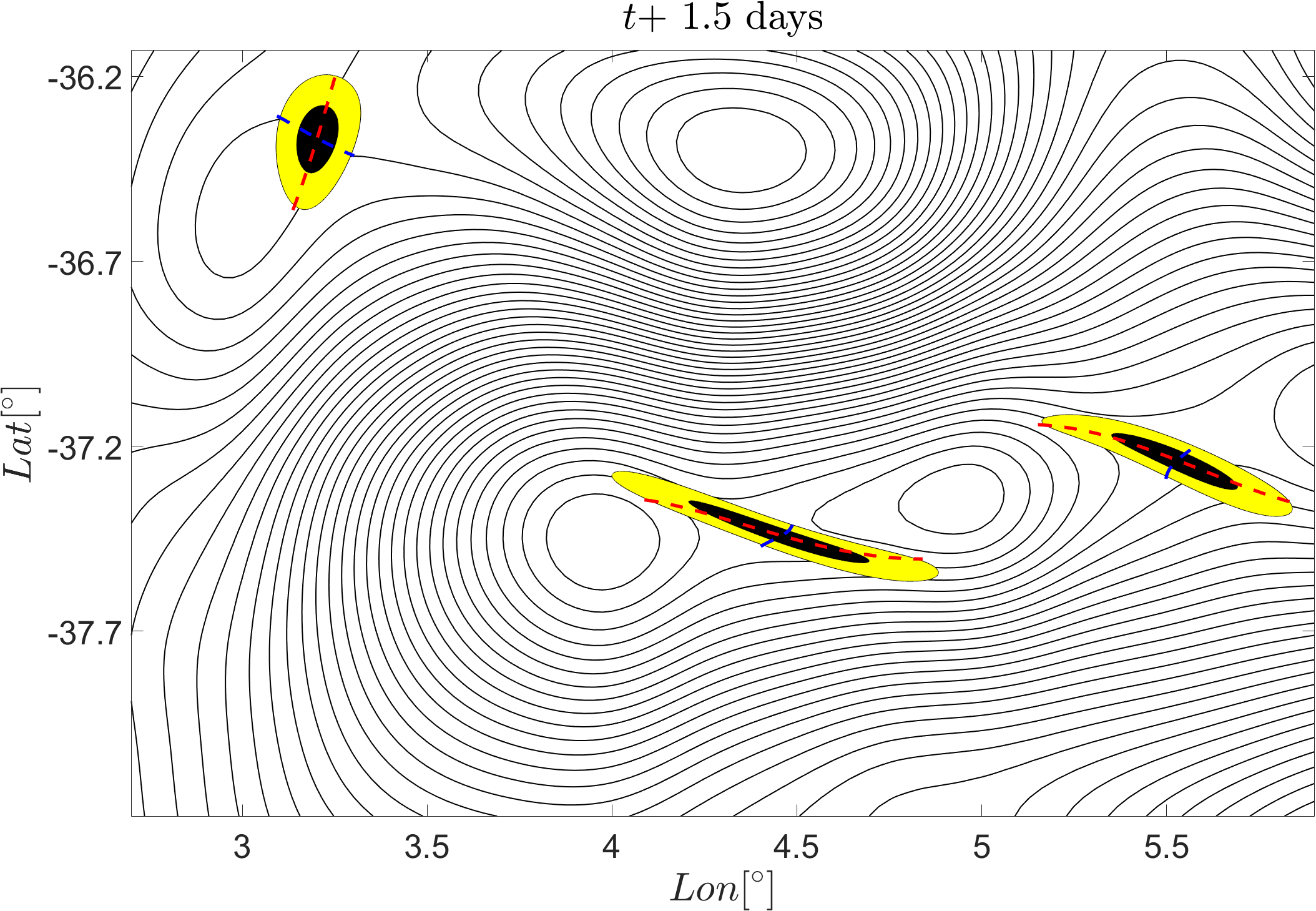}\label{fig:SaddHypStrml125}}\\
	\caption{(a) Saddle-type stagnation points of the velocity field (magenta triangles)
		with their associated stable (dashed blue) and unstable (dashed red)
		streamlines. Red dots denote the \textit{\emph{objective saddle points}}
		at the same time instant. (b-d) Advected images of the stable and
		unstable streamlines of saddle-type stagnation points up to 1.5 days
		and their effects on nearby particles.}

	\label{fig:CriticalPtsonStreamlines}
\end{figure}

A common Eulerian diagnostic for short-term hyperbolic tracer behavior
is the identification of instantaneous saddle-type stagnation points
in the velocity field. For comparison with this diagnostic, we show
in Fig. \ref{fig:SaddHypStrml0} the only three saddle-type stagnation
points (magenta triangles) that exist instantaneously in the domain
at time $t=24\ \mathrm{November}\ 2006$. Also shown are the corresponding
stable (dash blue) and unstable (dash red) directions inferred from
streamlines (black lines) for these stagnation points, along with
the ten \textit{\emph{objective saddle points}} (red dots) found at
the same time instant.

Two of these red dots marking hyperbolic OECS cores in Fig. \ref{fig:CriticalPtsonStreamlines}
fall near two instantaneous stagnation points, and give an improved
objective prediction for the cores of short-term saddle-type behavior
in Lagrangian particle motion. The improvement is seen by tracking
the deformation of the advected fluid blobs, initially centered on
the stagnation points. These blobs must stretch as they lie close
to two \textit{\emph{objective saddle points}}. The stretching blobs,
however, align more closely with the advected $e_{2}$-lines shown
in Fig. \ref{fig:Hypei12ECSs} when compared to the advected streamline
segments shown in Fig. \ref{fig:CriticalPtsonStreamlines}. The
	instantaneous stagnation point on the top left, instead, induces no
	notable material stretching on nearby particles since there are no
	hyperbolic OECS cores in its vicinity. The remaining eight hyperbolic
OECSs remain completely hidden in the instantaneous streamline picture
(cf. Fig. \ref{fig:CriticalPtsonStreamlines}), even though they induce
significant saddle-type material stretching, as we have already seen
in Fig. \ref{fig:Hypei12ECSs}.

These results illustrate that an instantaneous forecast strategy based
on saddle-type stagnation points of the velocity field may miss the
majority of significant short-term material stretching events in the
flow. Even for the detected stagnation points, there is no guarantee
that they signal nearby Lagrangian hyperbolic behavior correctly,
unless the velocity field has slow enough variation in their vicinity
\cite{HallerPoje1998}. This mismatch between objective and frame-dependent
predictions for saddle points is generally expected to increase further
for highly unsteady flows. \\

\subsection{Parabolic OECSs}

We now discuss the existence of parabolic OECSs in the full domain
used for computing elliptic OECSs. There is no known persistent Eulerian
or Lagrangian jet in this part of the ocean, but we nevertheless uncover
parabolic OECSs in this region that act as short-term pathways for
material transport.

\begin{figure}[h]
	\centering
\includegraphics[width=.8\textwidth]{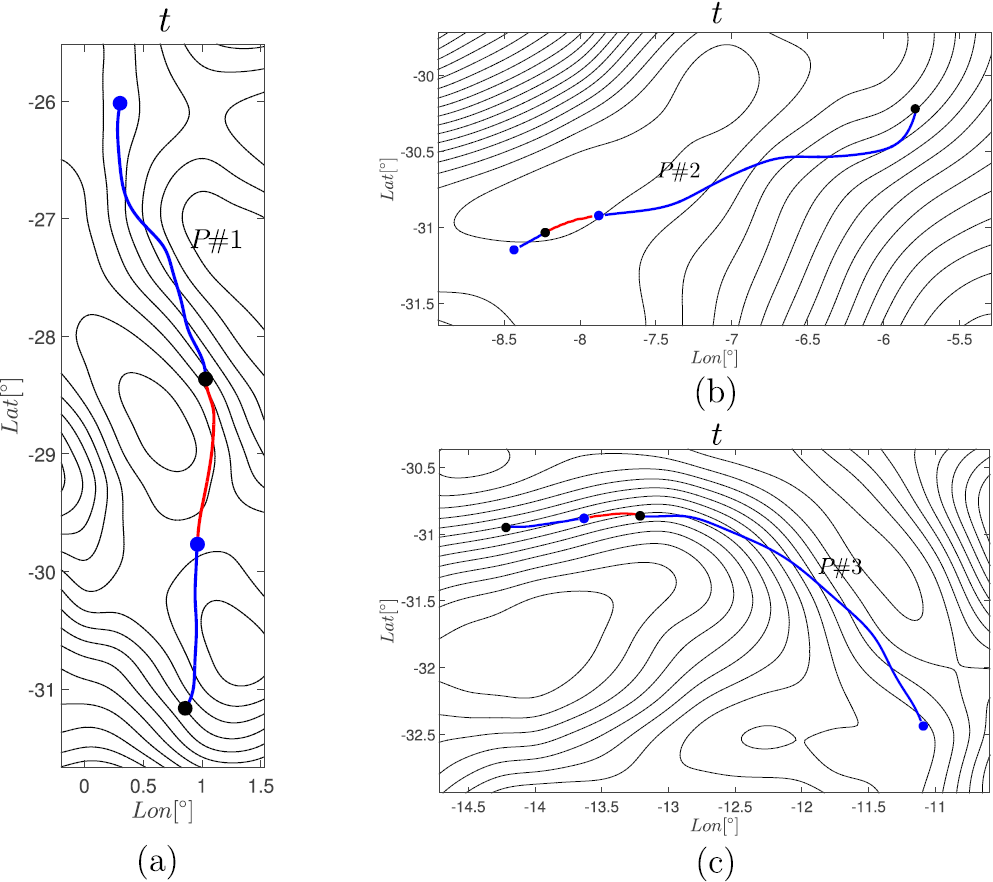}	
	\caption{(a-c) Parabolic OECSs as alternating heteroclinic connections of $e_{1}$-lines
		(blue) and $e_{2}$-lines (red) between trisector-type singularities
		(blue dots) and wedge-type singularities (black dots). Instantaneous
		streamlines are shown in the background.}
	\label{fig:ParabolicECSsOverview}
\end{figure}

In Fig. \ref{fig:ParabolicECSsOverview}, we show three parabolic
OECSs obtained from the application of Algorithm 4 to the velocity
data studied here. These OECSs are, therefore, alternating chains
of $e_{1}$-lines (blue) and $e_{2}$-lines (red), connecting trisector-type
(blue dots) singularities and wedge-type singularities (black dots),
with the streamlines shown in the background for reference. The three
parabolic OECSs are of approximate sizes $500$ km for $P\#1$ and
$300$ km for $P\#2$ and $P\#3$.

\begin{figure}[h]
	\subfloat[]{\includegraphics[width=0.45\textwidth]{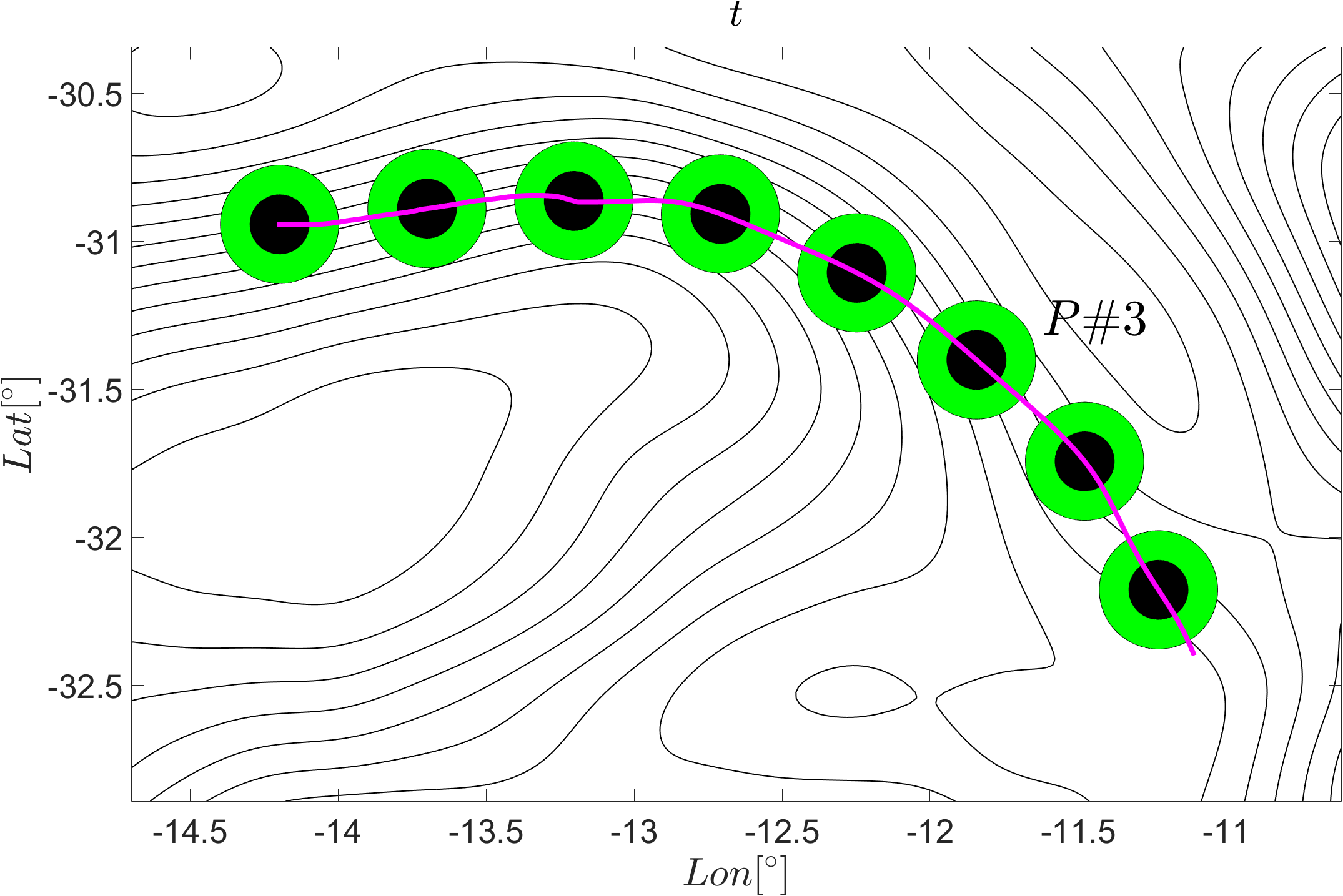}\label{fig:P3Strml0}}\hfill{}
	\subfloat[]{\includegraphics[width=0.45\textwidth]{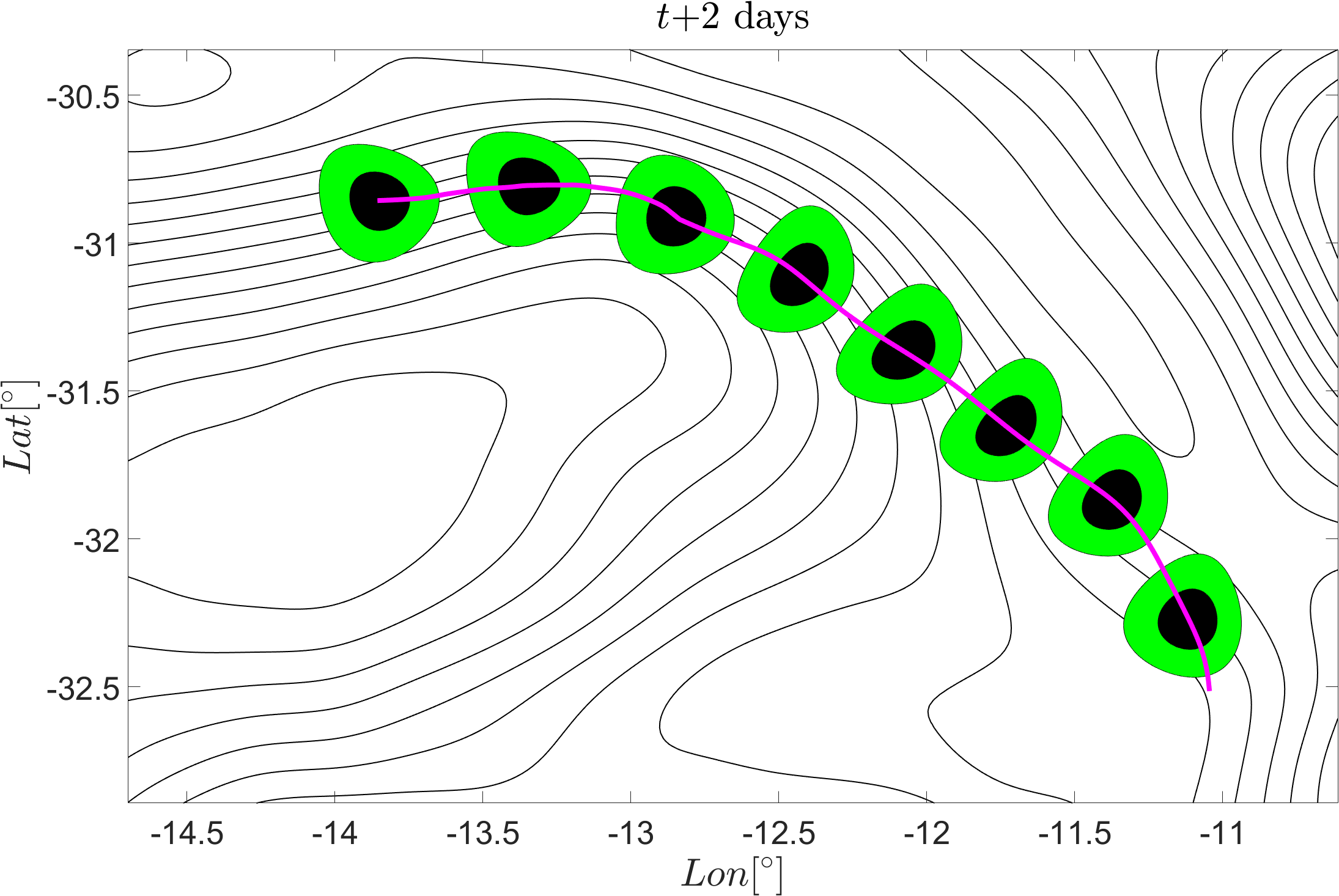}\label{fig:P3Strml2}}\\
	\subfloat[]{\includegraphics[width=0.45\textwidth]{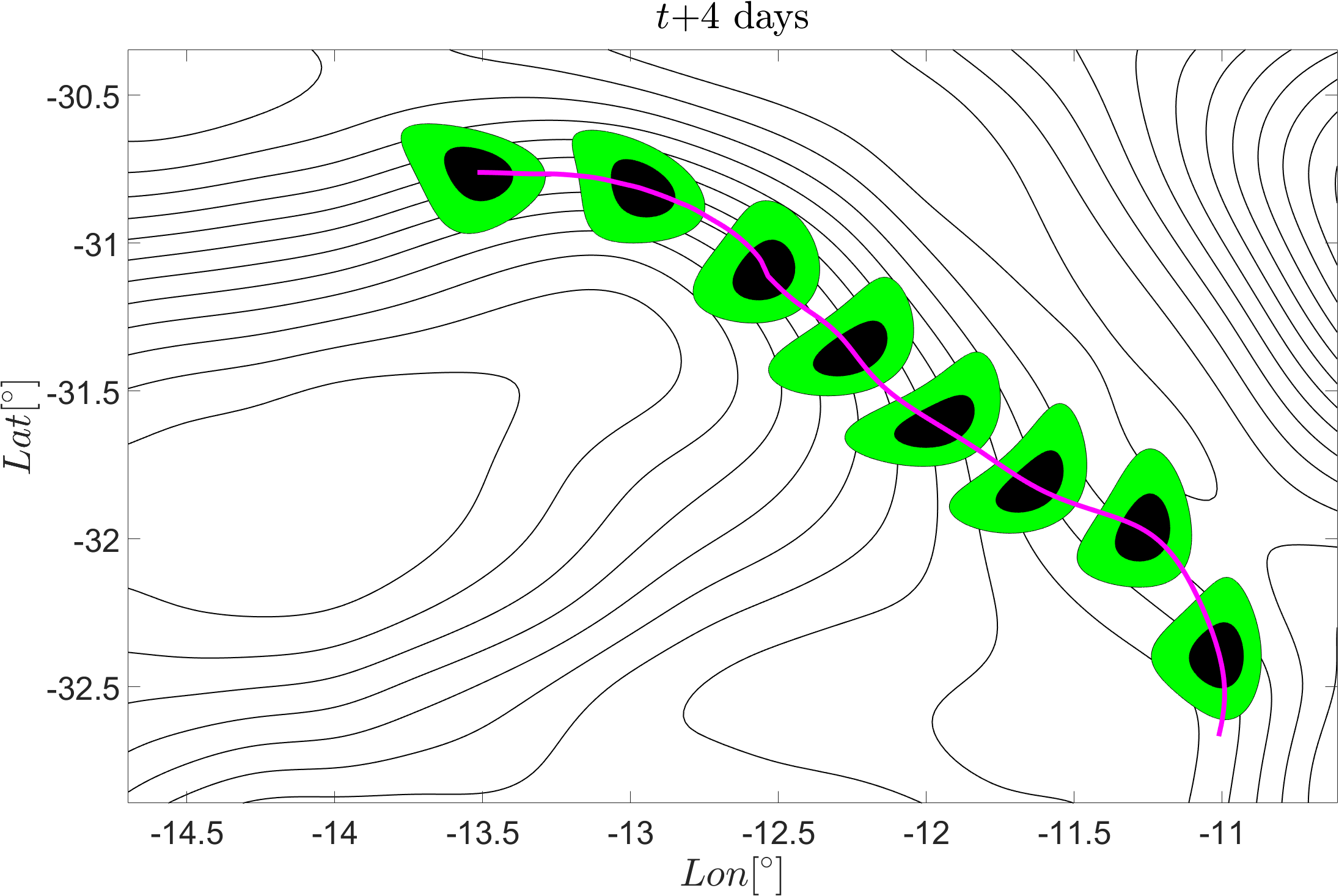}\label{fig:P3Strml4}}\hfill{}
	\subfloat[]{\includegraphics[width=0.45\textwidth]{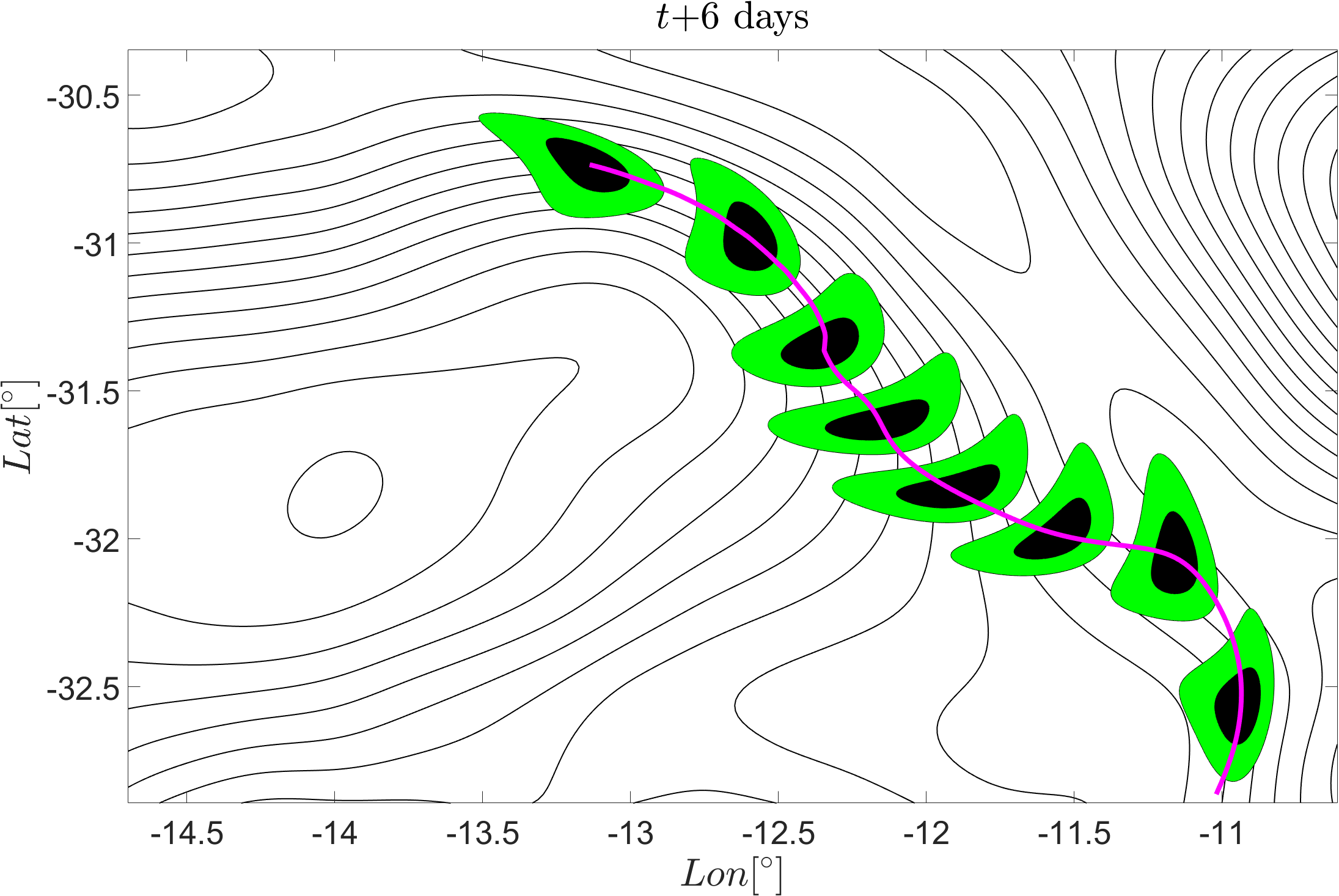}\label{fig:P3Strml6}}\\
	\caption{(a) Parabolic OECS, P\#3, at initial time
		overlaid on streamlines. (b-d) Advected images of the Parabolic OECS
		up to 6 days, and its effect on nearby particles.}
	\label{fig:ParabolicECS3adv}
\end{figure}

As a representative example, Fig. \ref{fig:ParabolicECS3adv} shows
the advected positions of the parabolic OECS, P\#3, along with the
deformation of material blobs initialized along this OECS, for integration
times up to six days. The advection illustrates the existence of a
short-term pathway along which initial conditions march towards the
bottom right. The deformation of the blobs has a characteristic boomerang
shape, which is similar to that observed along persistent Lagrangian
jets \cite{ShearlessBarrierFarazmand2014,JupiterAlirezaGeorge2014}.
This objective short-term transport route, therefore, has a clearly
verifiable material impact even though it has no clear signature in
the instantaneous streamline geometry.

\clearpage
\section{Conclusions}

We have developed a variational theory of objective Eulerian Coherent
Structures (OECSs) for two-dimensional, non-autonomous dynamical systems.
In this theory, we define OECSs at time $t$ as curves with the lowest
instantaneous material stretching rate or shearing rate in the flow. We
defined elliptic OECSs as closed stationary curves of the instantaneous
stretching-rate and shearless OECSs as stationary curves of the instantaneous
shearing-rate functional. A further classification of shearless OECSs
divides them into hyperbolic (attracting or repelling) and parabolic
(jet-type) OECSs. The approach we have taken here is objective, i.e.,
returns the same OECSs in frames translating and rotating relative
to each other. 

In our present, instantaneous context, objectivity guarantees a self-consistent
detection of short-term material coherent structures via OECSs in
unsteady, two-dimensional fluid flows. Indeed, we have shown that
elliptic OECS provide accurate identification of short-term material
vortices; hyperbolic OECSs reveal generalized saddle points with the
corresponding stable and unstable manifolds; and parabolic OECSs uncover
short-term jet-type pathways for material transport. We have verified
the accuracy of these OECS-based predictions by actual material advection
in our ocean data example obtained from satellite altimetry. 

We have also compared our results to other broadly used instantaneous
coherent structure indicators: streamline topology (which is neither
Galilean invariant nor objective) and the Okubo--Weiss criterion (which
is Galilean invariant but not objective). Such diagnostics often need
an ad hoc selection of threshold parameters, which limits the reliability
of the results they provide. In contrast, our procedure uses no free
parameters or thresholds.

We have found examples of false positives and false negatives suggested
by these non-objective Eulerian indicators. For instance, several
regions of maximum $OW$-- ellipticity turn out to stretch significantly
more than the regions identified by elliptic OECSs.\textbf{ }Even
when the indicators happen to suggest an OECS, the objective variational
principles developed here give more accurate results, as confirmed
by short-term material advection of the detected sets. 

As an alternative, elliptic OECSs can also be sought as closed curves
showing short-term rotational coherence, i.e., equal material rotation
rate \cite{HallerLAVD2015}. Obtained as the instantaneous limit of
a Lagrangian rotational coherence principle, rotationally coherent
OECSs are also objective. These OECSs do not restrict stretching rates
and hence would generally be expected to give slightly larger but
less coherent Eulerian vortices than the ones detected by the approach
developed here. Those larger vortices display tangential filamentation,
while the stretching-rate-based OECSs introduced here are instantaneous
limits of the perfectly coherent, black-hole-type elliptic LCSs derived
in \cite{BlackHoleHaller2013}. 

OECS-based forecasting of material transport and mixing is necessarily
confined to shorter time scales. Such shorter time scales, however,
are precisely the relevant ones for flow control or environmental
assessment where quick operational decisions need to be made. Results
on the use of OECSs in short-term forecasting and an extension to
higher dimensions will appear elsewhere.

\section*{Acknowledgment }

We would like to acknowledge Mohammad Farazmand for helpful discussions
on the subject of this paper.
\clearpage
\appendix

\section{Deformation-rate measures for material curves\label{sec:derEulMeas} }

\begin{figure}[h]
	\subfloat[]{\includegraphics[width=0.6\columnwidth]{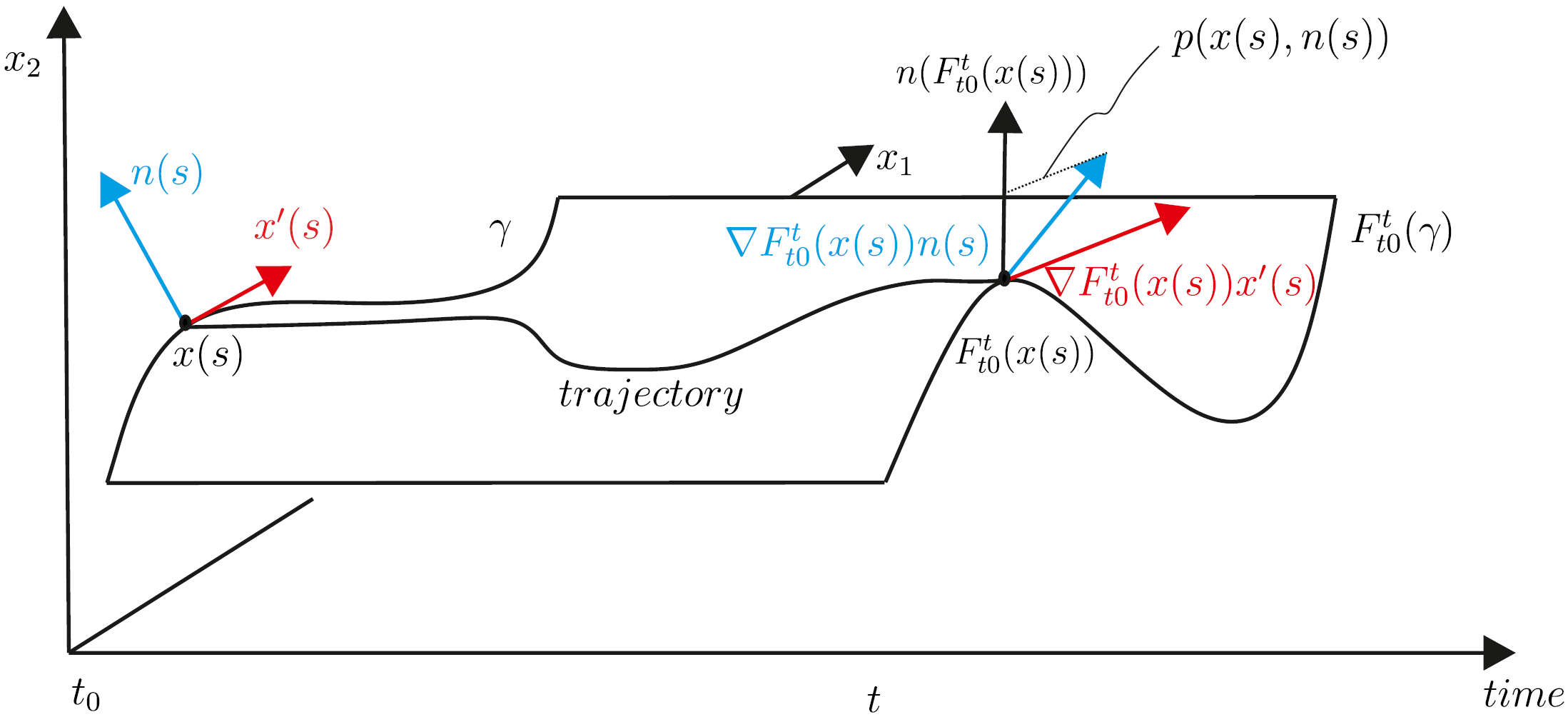}\label{fig:DefTpMLagr-1}}\hfill{}
	\subfloat[]{\includegraphics[width=0.3\columnwidth]{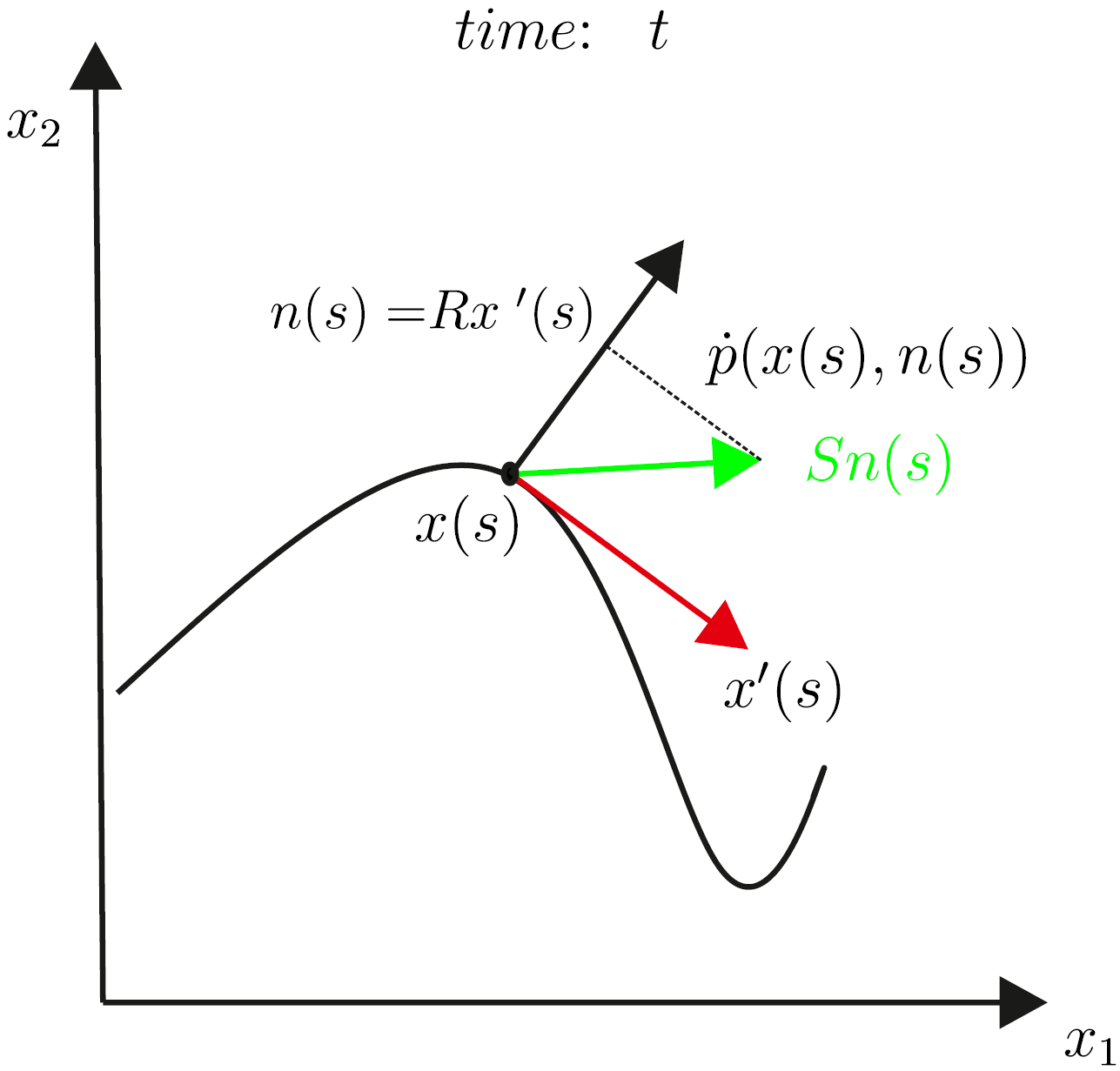}
		
		\label{fig:DefTpMEulr-1}}
	
	\caption{(a) Evolution of local tangent and normal vectors under the linearized
		flow map $\nabla F_{t_{0}}^{t}$ in the extended phase space of $x$
		and $t$. (b) Eulerian counterpart letting $x(s)$ be the unit speed
		parametrization of a regular curve at the current time $t$.}

	\label{fig:DefTgSp-1}
\end{figure}

At time $t_{0}$, consider a smooth curve of initial conditions $\gamma$,
parametrized as $s\mapsto x(s)$ via its arclength $s\in[0,\sigma]$.
Let $x'(s)$ and $n(s)=R\tfrac{x'(s)}{\vert x'(s)\vert}$ denote the
local tangent and normal vectors to $\gamma$ respectively. While
tangent vectors of $\gamma$ are mapped into tangent vectors by the
linearized flow map $\nabla F_{t_{0}}^{t}$, initial normal vectors
at time $t_{0}$ are not mapped by $\nabla F_{t_{0}}^{t}$ into the
normal space at time $t$ (Fig. \ref{fig:DefTpMLagr-1}). The local
deformation of $\gamma$ and its nearby trajectories, over the time
interval $[t_{0},t]$, can be expressed in terms of two Lagrangian
quantities: the tangential shear and the tangential strain over the
time interval $[t_{0},t]$ (see \cite{ShearlessBarrierFarazmand2014}
and \cite{BlackHoleHaller2013}).

Specifically, the tangential shear at point $x(s)$ is given by
\[
p_{t_{0}}^{t}(s)=D_{t_{0}}^{t}(x(s))=\dfrac{\langle x'(s),D_{t_{0}}^{t}(x(s))x'(s)\rangle}{\sqrt{\langle x'(s),C_{t_{0}}^{t}(x(s))x'(s)\rangle\langle x'(s),x'(s)\rangle}},\qquad\qquad D_{t_{0}}^{t}:=\dfrac{1}{2}[C_{t_{0}}^{t}R-RC_{t_{0}}^{t}],
\]
and the tangential strain at the same point is give by 
\[
q_{t_{0}}^{t}(s)=\dfrac{\sqrt{\langle x'(s),C_{t_{0}}^{t}(x(s))x'(s)\rangle}}{\sqrt{\langle x'(s),x'(s)\rangle}}.
\]

The quantities $p_{t_{0}}^{t}$ and $q_{t_{0}}^{t}$ give pointwise
information on tangential shear and tangential stretching experienced
by a material curve $\gamma$ over the time interval $[t_{0},t]$
in an objective (frame-invariant) fashion. 

The instantaneous rates of shear and tangential stretching along $\gamma$
can be obtained by differentiating $p_{t_{0}}^{t}$ and $q_{t_{0}}^{t}$
with respect to $t$ and setting $t=t_{0}$. Specifically, we have
\[
\dot{p}(s,t):=\tfrac{d}{d\tau}p_{t}^{t+\tau}(s)\vert_{\tau=0}=\dfrac{\langle x'(s),[S(x(s),t)R-RS(x(s),t)]x'(s)\rangle}{\langle x'(s),x'(s)\rangle},
\]
and 
\[
\begin{aligned}\dot{q}(s,t):=\tfrac{d}{d\tau}q_{t}^{t+\tau}(s)\vert_{\tau=0}= & \dfrac{1}{{\sqrt{\langle x'(s),x'(s)\rangle}}}\dfrac{\langle x'(s),\tfrac{d}{d\tau}C_{t}^{t+\tau}(x(s))\bigr\lvert_{\tau=0}x'(s)\rangle}{2\ \sqrt{{\langle x'(s),C_{t}^{t}(x(s))x'(s)\rangle}}}\\
= & \dfrac{\langle x'(s),S(x(s),t)x'(s)\rangle}{\langle x'(s),x'(s)\rangle},
\end{aligned}
\]
where we used the following relations 

\begin{equation}
\begin{aligned}\nabla F_{t}^{t}(x)= & I,\\
\tfrac{d}{d\tau}\nabla F_{t}^{t+\tau}(x)\bigr\lvert_{\tau=0}= & \nabla v(x,t)\nabla F_{t}^{t}(x)\\
= & \nabla v(x,t),\\
\tfrac{d}{d\tau}C_{t}^{t+\tau}(x)\bigr\lvert_{\tau=0}= & (\dot{\nabla F}_{t}^{t+\tau}(x)^{\top}\nabla F_{t}^{t+\tau}(x)+\nabla F_{t}^{t+\tau}(x)^{\top}\dot{\nabla F}_{t}^{t+\tau}(x))\bigr\lvert_{\tau=0}\\
= & 2S(x,t),\\
\tfrac{d}{d\tau}D_{t}^{t+\tau}\bigr\lvert_{\tau=0}(x)= & [S(x,t)R-RS(x,t)],\\
D_{t}^{t}= & 0.
\end{aligned}
\end{equation}

\section{$\chi_{\mu}^{\pm}$ is a one-parameter family of rotated vector fields\label{sec:RotVecFld} }

Assume there exists a limit cycle of (\ref{eq:chimufield}), $\gamma$,
for one choice of $\pm$ and a fixed value of $\mu$. The hyperbolic
nature of limit cycles guarantees their persistence with respect to
small changes in the parameter $\mu$, which leads to a one-parameter
family of limit cycles for the vector field $\chi_{\mu}^{\pm}$. In
general, these limit cycles can deform in an arbitrary fashion and
even intersect each other. 

In the present case, however, $\chi_{\mu}^{\pm}$ turns out to be
as a\textit{ one-parameter family of rotated vector fields }in the
sense of \cite{Duff1953}. This means that trajectories of (\ref{eq:chimufield}),
for each of the choices $\pm$, cannot intersect when $\mu$ varies
and shrink or expand for monotonic changes of the parameter. Hence,
limit cycles of \textit{one-parameter family of rotated vector fields,}
corresponding to different values of $\mu$, cannot intersect each
other.

To qualify as a\textit{ one-parameter family of rotated vector fields,
}$\chi_{\mu}^{\pm}$ must be locally smooth in a neighborhood of the
limit cycle, and the vector field defined by $\tfrac{d}{d\mu}\chi_{\mu}^{\pm}\times\chi_{\mu}^{\pm}$
must keep the same orientation with respect to the plane spanned by
$\tfrac{d}{d\mu}\chi_{\mu}^{\pm}$ and $\chi_{\mu}^{\pm}$.

Indeed, $\chi_{\mu}^{\pm}$ is smoothly orientable in the vicinity
of limit cycles, although it is not globally orientable due to orientational
discontinuities of the $e_{i}$ fields. For the second condition above,
it is equivalent to check that $\mathrm{sign}\langle\tfrac{d}{d\mu}\chi_{\mu}^{\pm}\times\chi_{\mu}^{\pm},e_{3}\rangle$,
with $e_{3}$ denoting the planar unit vector, remains unaltered over
the domain $U_{\mu}$ for each of the choices $\pm$. 

In the $[e_{1},e_{2},e_{3}]$ basis, $\langle\tfrac{d}{d\mu}\chi_{\mu}^{\pm}\times\chi_{\mu}^{\pm},e_{3}\rangle$
can be computed as 
\[
\begin{aligned}O_{\mu}^{\pm}(x)= & \chi_{\mu}^{\pm}\times\dfrac{d}{d\mu}\chi_{\mu}^{\pm}(x)=o_{\mu}^{\pm}(x)e_{3},\\
o_{\mu}^{\pm}(x)= & \langle\tfrac{d}{d\mu}\chi_{\mu}^{\pm}\times\chi_{\mu}^{\pm},e_{3}\rangle=\dfrac{\pm}{2\sqrt{(\mu-s_{1})(s_{2}-\mu)}},
\end{aligned}
\]
therefore, $\tfrac{d}{d\mu}\chi_{\mu}^{\pm}\times\chi_{\mu}^{\pm}$
keeps the same orientation for each of the signs $\pm$, under a monotonic
change of the parameter $\mu$.

With our choice of relative orientation between $e_{1}$ and $e_{2}$,
(cf. equation (\ref{eq:RelativeOrientation})), the sign of $o_{\mu}^{\pm}(x)$
gives the direction of rotation (positive counterclockwise) of the
$\chi_{\mu}^{\pm}$ field when the parameter $\mu$ increases, for
each of the choices $\pm$. Finally, the quantity $o_{\mu}^{\pm}(x)$
could be used in the computation of Elliptic OECSs through a systematic
change of the parameter $\mu$, as described in \cite{BlackHoleHaller2013}.

\newpage
\bibliographystyle{plain}

\bibliography{ReferenceList2}

\end{document}